\documentclass[11pt]{article}
\usepackage{verbatim}
\usepackage[dvips]{graphicx}
\usepackage{amsmath,amsthm,amssymb,amsbsy}
\usepackage{psfig,makeidx}
\usepackage{color}

\renewcommand{\baselinestretch}{1.5}

\newcommand{\bea}{\begin{eqnarray*}}
\newcommand{\eea}{\end{eqnarray*}}

\newcommand{\bydef}{\stackrel{\bigtriangleup}{=}}

\newcommand{\bc}{\begin{center}}
\newcommand{\ec}{\end{center}}

\newcommand{\R}{\mathbb R}

\newcommand{\be}{\begin{equation}}
\newcommand{\ee}{\end{equation}}
\newcommand{\ben}{\begin{enumerate}}
\newcommand{\een}{\end{enumerate}}
\newcommand{\db}{\hspace*{\fill}{\zapf o}}

\newcommand{\bpn}{\begin{proposition}\twlsf}
\newcommand{\epn}{\db\end{proposition}}
\newcommand{\bdm}{\begin{displaymath}}
\newcommand{\edm}{\end{displaymath}}

\newcommand{\mb}[1]{\boldsymbol{#1}}

\newcommand{\st}{\mathop{\rm s.t.}}

\newtheorem{proposition}{Proposition}

\newtheorem{theorem}{Theorem}

\theoremstyle{remark}

\theoremstyle{plain}

\theoremstyle{definition}

\DeclareMathAlphabet{\mathsfsl}{OT1}{cmss}{m}{sl}
\renewcommand{\v}[1]{\boldsymbol{#1}}

\newcommand{\field}[2]{\mathbb{#1}^{#2}}

\newcommand{\ch}[1]{\mathrm{conv}\left(#1\right)}

\newcommand{\trace}[1] {\mathrm{Tr}(#1)}
\newcommand{\prob}[1] {\mathbb{P}\left(#1\right)}
\newcommand{\ev}[2] {\mathbb{E}_{#1}\left[#2\right]}

\title{Theory and applications of Robust Optimization}
\author{Dimitris Bertsimas\footnote{Sloan School of Management and Operations Research Center,
Massachusetts Institute of Technology, E40-147, Cambridge, MA 02139.
\texttt{dbertsim@mit.edu}}, \ \ \ David B. Brown\footnote{Fuqua School of
Business, Duke University, 100 Fuqua Drive, Box 90120, Durham, NC 27708.
\texttt{dbbrown@duke.edu}}, \ \ \ Constantine
Caramanis\footnote{Department of Electrical and
Computer Engineering, The University of Texas at Austin, 1 University
Station, Austin, TX 78712.
\texttt{caramanis@mail.utexas.edu}}}
\date{First version: September 5, 2008\ \\ Revised version: \today}

\begin{document}
\maketitle

\bibliographystyle{plain}

\begin{abstract}
In this paper we survey the primary research, both theoretical and
applied, in the area of Robust Optimization (RO). Our focus is
on the computational attractiveness of RO approaches, as well as the modeling power and broad applicability of the methodology. In
addition to surveying prominent theoretical results of RO, we also present some recent results
linking RO to adaptable models for multi-stage decision-making
problems. Finally, we highlight applications of RO
across a wide spectrum of domains, including finance, statistics, learning, and various areas of engineering.\ 
\end{abstract}

\vspace{2.0in}

\thispagestyle{empty}

\pagebreak

\setcounter{page}{1}

\section{Introduction}\label{introsec}

Optimization affected by parameter uncertainty has long been a focus of the mathematical programming community. Solutions to optimization problems
can exhibit remarkable sensitivity to perturbations in the parameters of the problem (demonstrated in compelling fashion in \cite{BenTalNemirovski00}) thus often rendering a computed solution highly infeasible, suboptimal, or both (in short, potentially worthless).

In science and engineering, this is hardly a new notion. In the context of optimization, the most closely related field is that of Robust Control (we refer to the textbooks \cite{ZhouDoyleGlover1996} and \cite{DullerudPaganini1999}, and the references therein).
While there are many high-level similarities, and indeed much of the motivation for the development of Robust Optimization came from the Robust Control community, Robust Optimization is a distinct field, focusing on traditional optimization-theoretic concepts, particularly algorithms, geometry, and tractability, in addition to modeling power and structural results which are more generically prevalent in the setting of robustness.

Stochastic Optimization starts by assuming the uncertainty has a probabilistic description. This approach has a long and active history dating at least as far back as Dantzig's original paper \cite{Dantzig1955}. We refer the interested reader to several textbooks (\cite{Infanger1994, BirgeLouveaux97, Prekopa95,KallWallace94}) and the many references therein for a more comprehensive picture of Stochastic Optimization.

This paper considers Robust Optimization (RO), a more recent approach to optimization under uncertainty, in which the uncertainty model is not stochastic, but rather deterministic and set-based. Instead of seeking to immunize the solution in some probabilistic sense to stochastic uncertainty, here the decision-maker constructs a solution that is feasible for \emph{any} realization of the uncertainty in a given set. The motivation for this approach is twofold. First, the model of set-based uncertainty is interesting in its own right, and in many applications is an appropriate notion of parameter uncertainty. Second, computational tractability is also a primary motivation and goal. It is this latter objective that has largely influenced the theoretical trajectory of Robust Optimization, and, more recently, has been responsible for its burgeoning success in a broad variety of application areas.
%
The work of Ben-Tal and Nemirovski (e.g., \cite{BenTalNemirovski98,
BenTalNemirovski99, BenTalNemirovski00}) and El Ghaoui et al.
\cite{ElGhaouiLebret97,ElGhaouiOustryLebret98} in the late 1990s,
coupled with advances in computing technology and the
development of fast, interior point methods for convex optimization, particularly for semidefinite optimization (e.g., Boyd and
Vandenberghe, \cite{BoydVandenberghe96}) sparked a massive flurry of interest in the field of Robust Optimization.

Central issues we seek to address in this paper include tractability of robust optimization models; conservativeness of the RO formulation, and flexibility to apply the framework to different settings and applications.
We give a summary of the main issues raised, and results presented.

\begin{enumerate}
\item Tractability: In general, the robust version of a tractable\footnote{Throughout this paper, we use the term ``tractable'' repeatedly. We use this as shorthand to refer to problems that can be reformulated into equivalent problems for which there are known solution algorithms with worst-case running time polynomial in a properly defined input size (see, e.g., Section 6.6 of Ben-Tal and Nemirovski \cite{BentalElGhaouiNemirovskiBook}). Similarly, by ``intractable'' we mean the existence of such an algorithm for general instances of the problem would imply P=NP.}
optimization problem may not itself be tractable. We outline tractability results, which depend on the structure of the
nominal problem as well as the class of uncertainty set. Many
well-known classes of optimization problems, including LP, QCQP,
SOCP, SDP, and some discrete problems as well, have a RO formulation that is tractable. Some care must be taken in the choice of the uncertainty set to ensure that tractability is preserved.

\item Conservativeness and probability guarantees: RO constructs solutions that are deterministically immune to realizations of the uncertain parameters in certain sets. This approach may be the only reasonable alternative when the parameter uncertainty is not stochastic, or if distributional information is not readily available. But even if there is an underlying distribution, the tractability benefits of the Robust Optimization approach may make it more attractive than alternative approaches from Stochastic Optimization. In this case, we might ask for probabilistic guarantees for the robust solution that can be computed {\it a priori}, as a function of the structure and size of the uncertainty set. In the sequel, we show that there are several convenient, efficient, and well-motivated parameterizations of different classes of uncertainty sets, that provide a notion of a {\it budget of uncertainty}. This allows the designer a level of flexibility in choosing the tradeoff between robustness and performance, and also allows the ability to choose the corresponding level of probabilistic protection. In particular, a perhaps surprising implication is that while the robust optimization formulation is inherently max-min (i.e., worst-case), the solutions it produces need not be overly conservative, and in many cases are very similar to those produced by stochastic methods.

\item Flexibility: In Section \ref{sec:structuretractability}, we discuss a wide array of optimization classes, and also uncertainty sets, and consider the properties of the robust versions. In the final section of this paper, we illustrate the broad modeling power of Robust Optimization by presenting a wide variety of applications. We also give pointers to some surprising uses of robust optimization, particularly in statistics, where robust optimization is used as a tool to imbue the solution with desirable properties, like sparsity, stability or statistical consistency. \end{enumerate}

The overall aim of this paper is to outline the development and main aspects of Robust Optimization, with an emphasis on its flexibility and structure. While the paper is organized around some of the main themes of robust optimization research, we attempt throughout to compare with other methods, particularly stochastic optimization, thus providing guidance and some intuition on when the robust optimization avenue may be most appropriate, and ultimately successful.


We also refer the interested reader to the recent book of Ben-Tal, El Ghaoui and Nemirovski \cite{BentalElGhaouiNemirovskiBook}, which is an excellent reference on Robust Optimization that provides more detail on specific formulation and tractability issues. Our goal here is to provide a more condensed, higher level summary of key methodological results as well as a broad array of applications that use Robust Optimization.

\subsection*{A First Example}

To motivate RO and some of the modeling issues at hand, we begin with an example from portfolio selection. The example is a fairly standard one. We consider an investor who is attempting to allocate one unit of wealth among $n$ risky assets with random return $\tilde{\v{r}}$ and a risk-free asset (cash) with known return $r_f$. The investor may not short-sell risky assets or borrow. His goal is to optimally trade off between expected return and the probability that his portfolio loses money.

If the returns are stochastic with known distribution, the tradeoff between expected return and loss probability is a stochastic program. However, calculating a point on the pareto frontier is in general NP-hard even when the distribution of returns is discrete (Benati and Rizzi \cite{BenatiRizzi}).

We will consider two different cases: one where the distribution of asset price fluctuation matches the empirical distribution of given historical data and hence is known exactly, and then the case where it only approximately matches historical data. The latter case is of considerable practical importance, as the distribution of new returns (after an allocation decision) often deviate significantly from past samples. We compare the stochastic solution to several easily solved RO-based approximations in both of these cases.

The intractability of the stochastic problem arises because of the probability constraint on the loss:
\begin{equation}
\label{eq:example-probconstraint}
\mathbb{P}(\tilde{\v{r}}'\v{x}+r_f(1-\mathbf{1}'\v{x})\geq 1) \geq 1 - p_{{\tiny{loss}}},
\end{equation}
where $\v{x}$ is the vector of allocations into the $n$ risky assets (the decision variables). The robust optimization formulations replace this probabilistic constraint with a {\it deterministic constraint}, requiring the return to be nonnegative {\it for any realization of the returns in some given set}, called the uncertainty set:
\begin{equation}
\label{eq:example-robconstraint}
\tilde{\v{r}}'\v{x}+r_f(1-\mathbf{1}'\v{x})\geq 1 \ \forall\tilde{\v{r}}\in\mathcal{R}.
\end{equation}
While not explicitly specified in the robust constraint (\ref{eq:example-robconstraint}), the resulting solution has some $p_{{\tiny{loss}}}$.
As a rough rule, the bigger the set $\mathcal{R}$, the lower the objective function (since there are more constraints to satisfy), and the smaller the loss probability $p_{{\tiny{loss}}}$. Central themes in robust optimization are understanding how to structure the uncertainty set $\mathcal{R}$ so that the resulting problem is tractable and favorably trades off expected return with loss probability $p_{{\tiny{loss}}}$. Section \ref{sec:structuretractability} is devoted to the tractability of different types of uncertainty sets. Section \ref{sec:probability} focuses on obtaining {\it a priori} probabilistic guarantees given different uncertainty sets.
%
%
%
%
%
%
Here, we consider three types of uncertainty sets, all defined with a parameter to control ``size'' so that we can explore the resulting tradeoff of return, and $p_{{\tiny{loss}}}$:
\begin{eqnarray*}
\mathcal{R}^Q(\gamma) &=& \left\{\tilde{\v{r}} \ : \ (\tilde{\v{r}}-\hat{\v{r}})'\v{\Sigma}^{-1}(\tilde{\v{r}}-\hat{\v{r}})\leq\gamma^2\right\}, \\
\mathcal{R}^D(\Gamma) &=& \left\{\tilde{\v{r}} \ : \exists\v{u}\in\field{R}{n}_+ \ \text{s.t.} \ \tilde{r}_i=\hat{r}_i+(\underline{r}_i-\hat{r}_i)u_i, \ u_i\leq 1, \sum\limits_{i=1}^nu_i\leq\Gamma\right\}, \\
\mathcal{R}^T(\alpha) &=& \left\{\tilde{\v{r}} \ : \ \exists \v{q}\in\field{R}{N}_+ \ \text{s.t.} \ \tilde{\v{r}}=\sum\limits_{i=1}^Nq_i\v{r}^i, \ \mathbf{1}'\v{q}=1, \ q_i\leq \frac{1}{N(1-\alpha)}, \ i=1,\ldots,N\right\}.
\end{eqnarray*}

The set $\mathcal{R}^Q(\gamma)$ is a quadratic or ellipsoidal uncertainty set: this set considers all returns within a radius of $\gamma$ from the mean return vector, where the ellipsoid is tilted by the covariance. When $\gamma=0$, this set is just the singleton $\{\hat{\v{r}}\}$. The set $\mathcal{R}^D(\Gamma)$ ($D$ for ``D-norm'' model considered in Section \ref{sec:structuretractability}) considers all returns such that each component of the return is in the interval $[\underline{r}_i,\hat{r}_i]$, with the restriction that the total weight of deviation from $\hat{r}_i$, summed across all assets, may be no more than $\Gamma$. When $\Gamma=0$, this set is the singleton $\{\hat{\v{r}}\}$; at the other extreme, when $\Gamma=n$, returns in the range $[\underline{r}_i,\hat{r}_i]$ for all assets are considered. Finally, $\mathcal{R}^T(k)$ is the ``tail'' uncertainty set, and considers the convex hull of all possible $N(1-\alpha)$ point averages of the $N$ returns. When $\alpha=0$, this set is the singleton $\{\hat{\v{r}}\}$. When $\alpha=(N-1)/N$, this set is the convex hull of all $N$ returns.


To illustrate the use of these formulations, consider $n=10$ risky assets based on $N=300$ past market returns. The assets are a collection of equity and debt indices, and the return observations are monthly from a data set starting in 1981. For each of the three uncertainty RO formulations, we solve $200$ problems, each maximizing expected return subject to feasibility and the robust constraint at one of $200$ different values of their defining parameter $\gamma$, $\Gamma$, or $\alpha$. In total, we solve $600$ RO formulations. For comparison, we also formulate the problem of minimizing probability of loss subject to an expected return constraint as a stochastic program (which can be formulated as a mixed integer program), and solve $8$ versions of this problem, each corresponding to one of $8$ different expected return levels. The computations are performed using the MOSEK optimization toolbox in Matlab on a laptop computer with a 2.13GHZ processor and 2GB of RAM.

\begin{figure}[!t]
\begin{center}$
\begin{array}{cc}
\hspace{-.25in}
\includegraphics[trim = 1.25in 3.4in 1.35in 3.5in, clip, scale=.4]{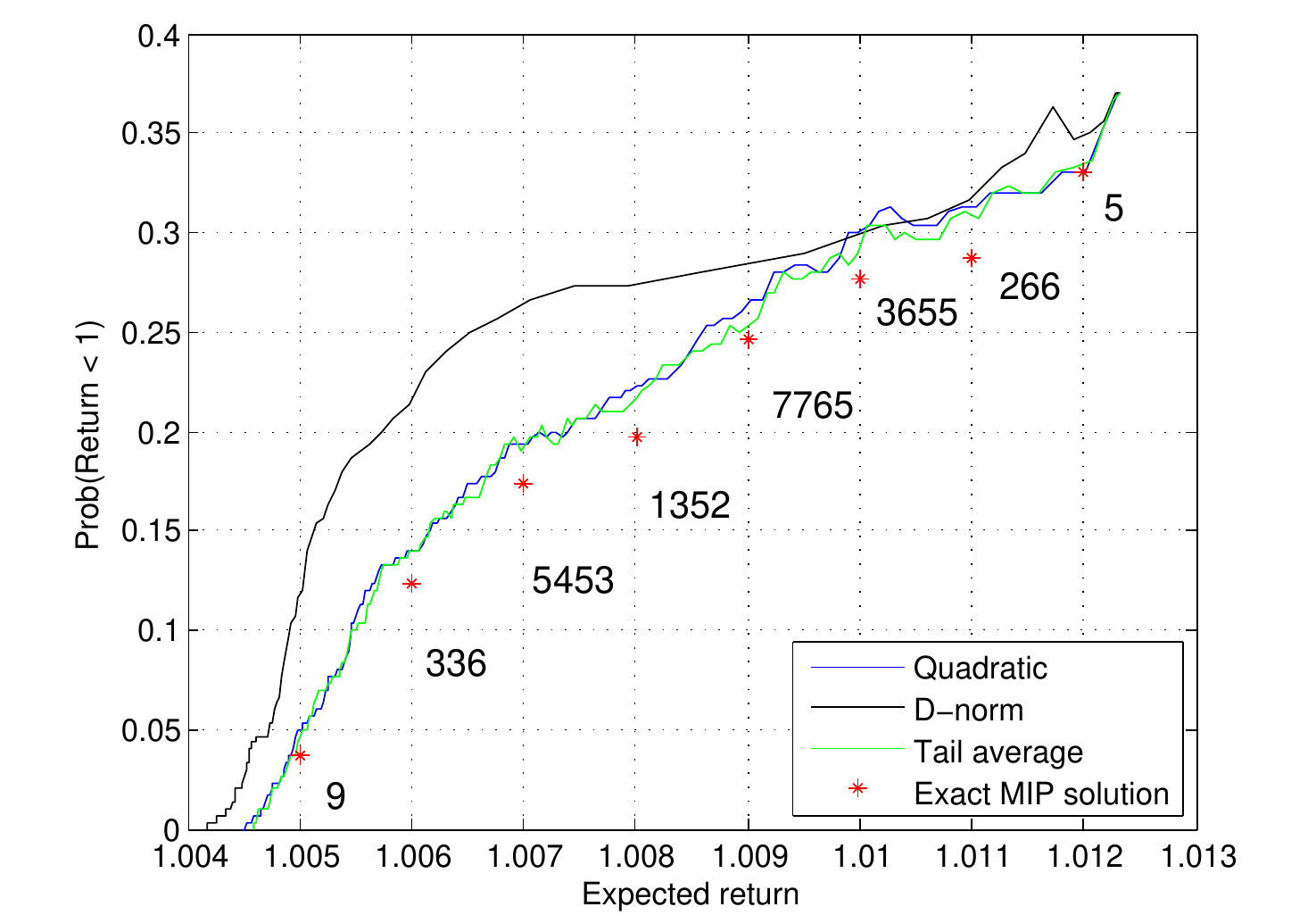} &
\includegraphics[trim = 1.25in 3.4in 1.35in 3.5in, clip, scale=.4]{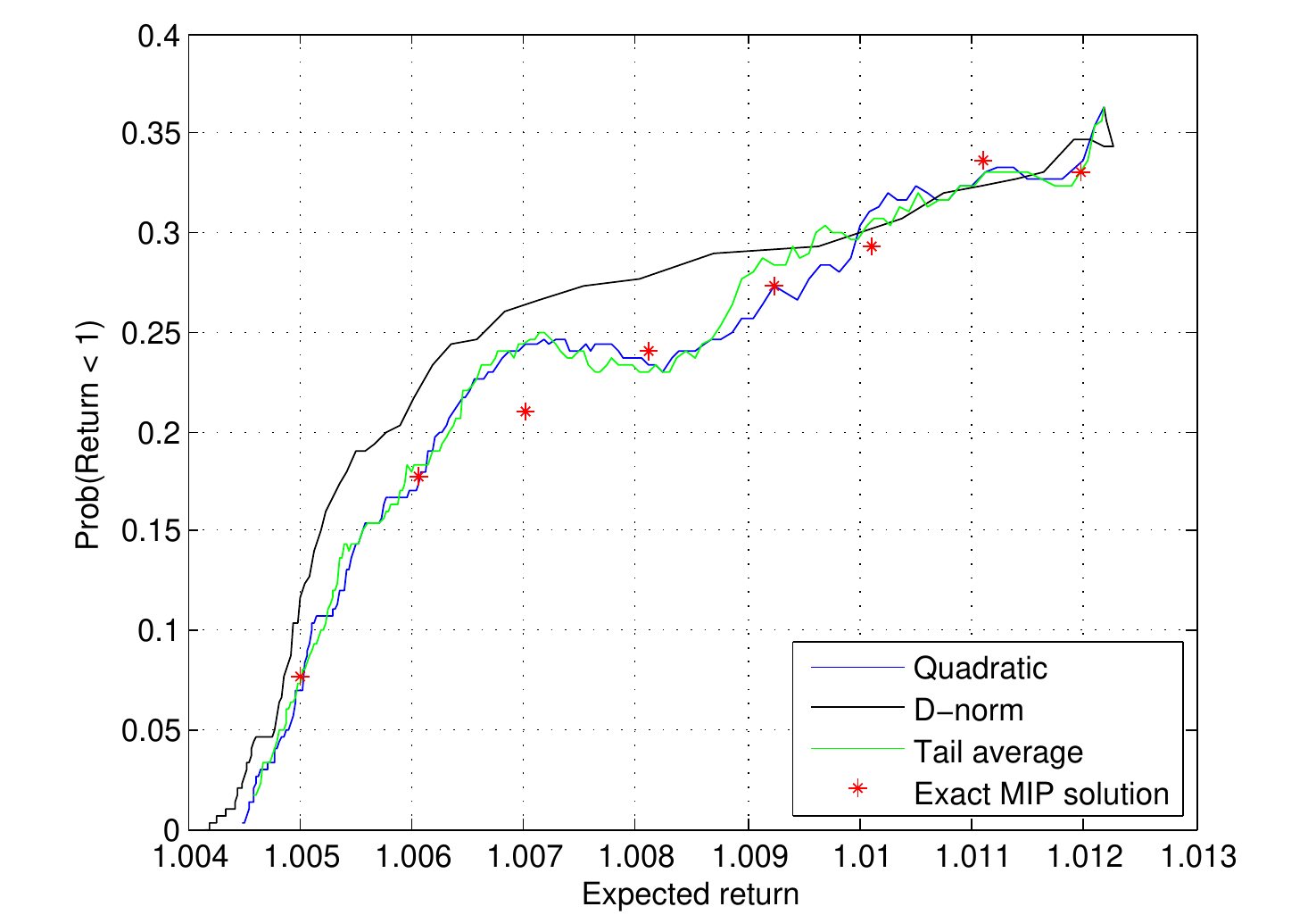}
\includegraphics[trim = 1.25in 3.4in 1.35in 3.5in, clip, scale=.4]{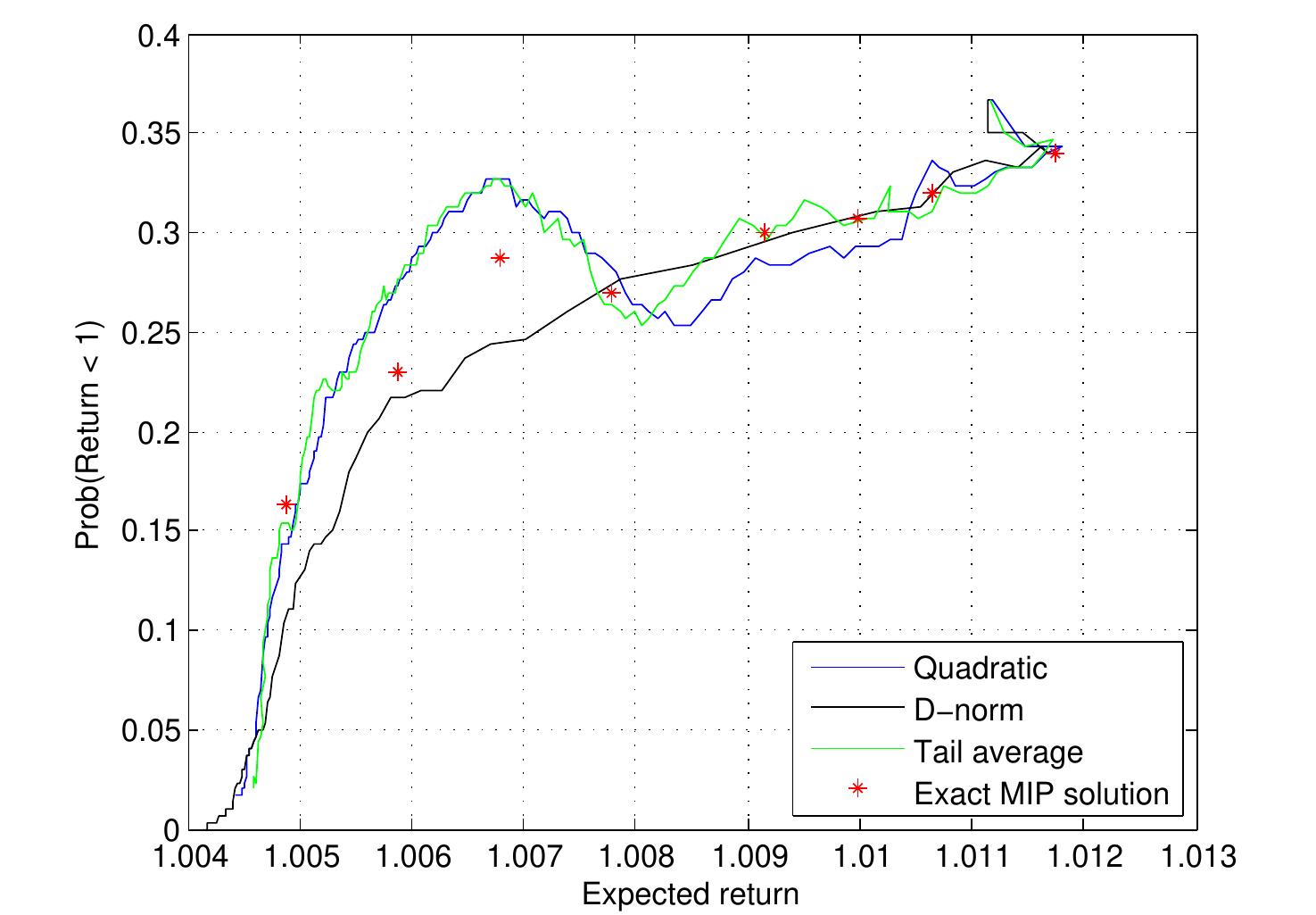}
\end{array}$
\end{center}
\vspace{-.25in}
\caption{(L): Expected return-loss probability frontier for RO-based formulations and exact stochastic formulation; numbers are time (sec.) for solving each stochastic program. (C/R): Frontier for model with random perturbations bounded by $1\%$ (C) and $2\%$ (R).}\label{examplefig}
\end{figure}

The results are shown in Figure \ref{examplefig}. On the left, we see the frontier for the three RO-based formulations as well as the performance of the exact formulation (at the $8$ return levels). The numbers indicate the time in seconds to solve the stochastic program in each case.

The stochastic model is designed for the nominal case, so we expect it to outperform the three RO-based formulations. However, even under this model, the gap from the $\mathcal{R}^Q$ and $\mathcal{R}^T$ RO frontiers is small: in several of the cases, the difference in performance is almost negligible. The largest improvement offered by the stochastic formulation is around a $2-3\%$ decrease in loss probability. Here, the solutions from the $\mathcal{R}^D$ model do not fare as well; though there is a range in which its performance is comparable to the other two RO-based models, typically its allocations appear to be conservative.
In general, solving the stochastic formulation exactly is difficult, which is not surprising given its NP-hardness. Though a few of the instances at extreme return levels are solved in only a few seconds, several of the instances require well over an hour to solve, and the worst case requires over $2.1$ hours to solve. The total time to solve these 8 instances is about $5.2$ hours; by contrast, solving the $600$ RO-based instances takes a bit under $10$ minutes in total, or about one second per instance.

On the center and right parts of Figure \ref{examplefig} are results for the computed portfolios under the same return model but with random perturbations. Specifically, we perturb each of the $N\times n$ return values by a random number uniformly distributed on $[.99,1.01]$ in the middle figure and $[.98,1.02]$ in the right figure. 
At the $1\%$ perturbation level, the gap in performance between the models is reduced, and there are regions in which each of the models are best as well as worst. The model based on $\mathcal{R}^D$ is least affected by the perturbation; its frontier is essentially unchanged. The models based on $\mathcal{R}^Q$ and $\mathcal{R}^T$ are more significantly affected, perhaps with the effect on $\mathcal{R}^T$ being a bit more pronounced. Finally, the stochastic formulation's solutions are the most sensitive of the bunch: though the SP solution is a winner in one of the $8$ cases, it is worse off than the others in several of the other cases, and the increase in loss probability from the original model is as large as $5-6\%$ for the SP solutions.

At the $2\%$ level, the results are even more pronounced: here, the SP solutions are always outperformed by one of the robust approaches, and the solutions based on $\mathcal{R}^D$ are relatively unaffected by the noise. The other two robust approaches are substantially affected, but nonetheless still win out in some parts of the frontier. When noise is introduced, it does not appear that the exact solutions confer much of an advantage and, in fact, may perform considerably worse. Though this is only one random trial, such results are typical.

There are several points of discussion here. First is the issue of complexity. The RO-based models are all fairly easy to solve here, though they themselves have complexities that scale differently. The $\mathcal{R}^Q$ model may be formulated as a second-order cone program (SOCP); both the $\mathcal{R}^D$ and the $\mathcal{R}^T$ models may be formulated as an LP. Meanwhile, the exact stochastic model is an NP-hard mixed integer program. Under the original model, it is clearly much easier to solve these RO-based models than the exact formulation. In a problem with financial data, it is easy to imagine having thousands of return samples. Whereas the RO formulations can still be solved quickly in such cases, solving the exact SP could be hopeless.

A second issue is the ability of solution methods to cope with deviations in the underlying model (or ``model uncertainty''). The RO-based formulations themselves are different in this regard. Here, the $\mathcal{R}^D$ approach focuses on the worst-case returns on a subset of the assets, the $\mathcal{R}^Q$ approach focuses on the first two moments of the returns, and the $\mathcal{R}^T$ approach focuses on averages over the lower tail of the distribution. Though all of these are somehow ``robust,'' $\mathcal{R}^D$ is the ``most robust'' of the three; indeed, we also implemented perturbations at the $5\%$ level and found its frontier is relatively unchanged, while the other three frontiers are severely distorted. Intuitively, we would expect models that are more robust will fare better in situations with new or altered data; indeed, we will later touch upon some work that shows that there are intimate connections between the robustness of a model and its ability to generalize in a statistical learning sense.

This idea - that Robust Optimization is useful in dealing with erroneous or noise-corrupted data - seems relatively well understood by the optimization community (those who build, study, and solve optimization models) at-large. In fact, we would guess that many figure this to be the \emph{raison d'\^{e}tre} for Robust Optimization. The final point that we would like to make is that, while dealing with perturbations is one virtue of the approach, RO is also more broadly of use as a computationally viable way to handle uncertainty in models that are \emph{on their own} quite difficult to solve, as illustrated here.

In this example, even if we are absolutely set on the original model, it is hard to solve exactly. Nonetheless, two of the RO-based approaches perform well and are not far from optimal under the nominal model. In addition, they may be computed orders of magnitude faster than the exact solution. Of course, we also see that the user needs to have some understanding of the structure of the uncertainty set in order to intelligently use RO techniques: the approach with $\mathcal{R}^D$, though somewhat conservative in the original model, is quite resistant to perturbations of the model. 

In short, RO provides a set of tools that may be useful in dealing with different types of uncertainties  - both the ``model error'' or ``noisy data'' type as well as complex, stochastic descriptions of uncertainty in an explicit model - in a computationally manageable way. Like any approach, however, there are tradeoffs, both in terms of performance issues and in terms of problem complexity. Understanding and managing these tradeoffs requires expertise. The goal of this paper, first and foremost, is to describe some of this landscape for RO. This includes detailing what types of RO formulations may be efficiently solved in large scale, as well what connections various RO formulations have to perhaps more widely known methods. The second goal of this paper is to then highlight an array of application domains for which some of these techniques have been useful.

\section{Structure and tractability results}
\label{sec:structuretractability} In this section, we outline
several of the structural properties, and detail some tractability
results of Robust Optimization. We also show how the notion of a
budget of uncertainty enters into several different uncertainty set
formulations.

\subsection{Robust Optimization}\label{rossec}
Given an objective $f_0(\v{x})$ to optimize, subject to constraints $f_i(\v{x},\v{u}_i) \leq 0$ with uncertain parameters, $\{\v{u}_i\}$,
%
the general Robust Optimization formulation is:
\begin{eqnarray}\label{rogenprob}
\text{minimize} && f_0(\v{x}) \nonumber \\
\text{subject to} && f_i(\v{x},\v{u}_i)\leq 0, \qquad \forall \
\v{u}_i\in\mathcal{U}_i, \ i=1,\ldots,m.
\end{eqnarray}
Here $\v{x}\in\field{R}{n}$ is a vector of decision variables,
$f_0,f_i: \mathbb{R}^n \rightarrow \mathbb{R}$ are functions, and the uncertainty parameters $\v{u}_i\in\field{R}{k}$ are assumed to take arbitrary values in the {\it uncertainty sets} $\mathcal{U}_i\subseteq\field{R}{k}$,
which, for our purposes, will always be closed. The goal of
(\ref{rogenprob}) is to compute minimum cost solutions $\v{x}^*$
among all those solutions which are feasible for \emph{all}
realizations of the disturbances $\v{u}_i$ within $\mathcal{U}_i$.
Thus, if some of the $\mathcal{U}_i$ are continuous sets,
(\ref{rogenprob}), as stated, has an infinite number of constraints. Intuitively, this problem offers some measure of feasibility
protection for optimization problems containing parameters which are not known exactly.

It is worthwhile to notice the following, straightforward facts
about the problem statement of (\ref{rogenprob}):
\begin{itemize}
\item The fact that the objective function is unaffected by
parameter uncertainty is without loss of generality; we may always
introduce an auxiliary variable, call it $t$, and minimize $t$
subject to the additional constraint
$\max\limits_{\v{u_0}\in\mathcal{U}_0}f_0(\v{x},\v{u}_0)\leq t$.

\item It is also without loss of generality to assume that the
uncertainty set $\mathcal{U}$ has the form
$\mathcal{U}=\mathcal{U}_1\times\ldots\times\mathcal{U}_m$. If we have a single uncertainty set $\mathcal{U}$ for which we
require $(\v{u}_1,\ldots,\v{u}_m)\in\mathcal{U}$, then the
constraint-wise feasibility requirement implies an equivalent
problem is (\ref{rogenprob}) with the $\mathcal{U}_i$ taken as the
projection of $\mathcal{U}$ along the corresponding dimensions (see
Ben-Tal and Nemirovski, \cite{BenTalNemirovski99} for more on this).

\item Constraints without uncertainty are also captured in this framework by
assuming the corresponding $\mathcal{U}_i$ to be singletons.

\item Problem (\ref{rogenprob}) also contains the instances when
the decision or disturbance vectors are contained in more general
vector spaces than $\field{R}{n}$ or $\field{R}{k}$ (e.g.,
$\field{S}{n}$ in the case of semidefinite optimization) with the
definitions modified accordingly.
\end{itemize}

Robust Optimization is distinctly different than \emph{sensitivity analysis}, which is typically applied as a post-optimization tool for quantifying the change in cost of the associated optimal solution with small perturbations in the underlying problem data. Here, our goal is to compute fixed solutions that ensure feasibility \emph{independent of the data}. In other words, such solutions have \emph{a priori} ensured feasibility when the problem parameters vary within the prescribed uncertainty set, which may be large. We refer the reader to some of the standard optimization literature (e.g., Bertsimas and Tsitsiklis, \cite{BertsimasTsitsiklis1997}, Boyd and Vandenberghe,
\cite{BoydVandenberghe2004}) and works on perturbation theory (e.g., Freund, \cite{Freund85}, Renegar, \cite{Renegar94}) for more on
sensitivity analysis.

It is not at all clear when (\ref{rogenprob}) is efficiently
solvable. One might imagine that the addition of robustness to a
general optimization problem comes at the expense of significantly
increased computational complexity. It turns out that this is true: the robust counterpart to an arbitrary convex optimization problem is in general intractable (\cite{BenTalNemirovski98}; some approximation results for robust convex problems with a conic structure are discussed in \cite{BertsimasSim05}).
Despite this, there are many robust problems that may be handled
in a tractable manner, and much of the literature has focused on specifying classes of functions $f_i$, coupled with the types of uncertainty sets $\mathcal{U}_i$, that yield tractable robust counterparts. If we define the robust feasible set to be
\begin{eqnarray}\label{rofeasset}
X(\mathcal{U}) &=& \left\{\v{x} \ | \ f_i(\v{x},\v{u}_i)\leq 0 \
\forall \ \v{u}_i\in\mathcal{U}_i, \ i=1,\ldots,m\right\},
\end{eqnarray}
then for the most part,\footnote{i.e., subject to a Slater condition.} tractability is tantamount to
$X(\mathcal{U})$ being convex in $\v{x}$, with an efficiently
computable separation test. More precisely, in the next section we
show that this is related to the structure of a particular
subproblem. We now present an abridged taxonomy of some of the main
results related to this issue.

\subsection{\textbf{Robust linear optimization}}\label{rlpsssec}

The robust counterpart of a linear optimization problem is written,
without loss of generality, as
\begin{eqnarray}\label{roblpbg}
\text{minimize} && \v{c}^{\top}\v{x} \nonumber \\
\text{subject to} && \v{Ax}\leq\v{b}, \qquad \forall \
\v{a_1}\in\mathcal{U}_1, \dots, \v{a_m}\in\mathcal{U}_m,
\end{eqnarray}
where $\v{a}_i$ represents the $i^{th}$ row of the uncertain matrix
$\v{A}$, and takes values in the uncertainty set
$\mathcal{U}_i\subseteq\field{R}{n}$. Then, $\mb{a}_i^{\top} \mb{x} \leq b_i$, $\forall \mb{a}_i
\in \mathcal{U}_i$, if and only if
$\max_{\{\mb{a}_i \in \mathcal{U}_i\}}
\mb{a}_i^{\top} \mb{x} \leq b_i$, $\forall \, i$.
We refer to this as the {\it subproblem} which must be solved. Ben-Tal and Nemirovski \cite{BenTalNemirovski99} show that the robust LP is essentially always tractable for most practical uncertainty sets of interest. Of course, the resulting robust problem may no longer be an LP. We now provide some more detailed examples.

{\bf Ellipsoidal Uncertainty:} Ben-Tal and Nemirovski
\cite{BenTalNemirovski99}, as well as El Ghaoui et al.
\cite{ElGhaouiLebret97,ElGhaouiOustryLebret98}, consider ellipsoidal
uncertainty sets. Controlling the size of these ellipsoidal sets, as in the theorem below, has the interpretation of a budget of uncertainty that the decision-maker selects in order to easily trade off robustness and performance.

\begin{theorem}\label{btlpthm}
(Ben-Tal and Nemirovski, \cite{BenTalNemirovski99}) Let
$\mathcal{U}$ be ``ellipsoidal,'' i.e.,
$$
\mathcal{U} = U(\Pi,\v{Q}) = \left\{\Pi(\v{u}) \ | \ \|\v{Qu}\|\leq \rho \right\},
$$
where $\v{u}\rightarrow\Pi(\v{u})$ is an affine embedding of
$\field{R}{L}$ into $\field{R}{m\times n}$ and
$\v{Q}\in\field{R}{M\times L}$.
Then Problem (\ref{roblpbg}) is equivalent to a second-order cone
program (SOCP). Explicitly, if we have the uncertain optimization
$$
\begin{array}{rl}
\mbox{minimize} & \mb{c}^{\top}\mb{x} \\
\mbox{subject to} & \mb{a}_i \mb{x} \leq b_i, \quad \forall \mb{a}_i \in {\cal U}_i, \quad \forall i=1,\dots,m,
\end{array}
$$
where the uncertainty set is given as:
$$
{\cal U} = \{(\mb{a}_1,\dots,\mb{a}_m) \,:\, \mb{a}_i = \mb{a}_i^0 + \Delta_i u_i, \,\,
i=1,\dots,m, \quad ||u||_2 \leq \rho\},
$$
($\mb{a}_i^0$ denotes the nominal value)
then the robust counterpart is:
\begin{eqnarray*}
\mbox{mininize} && \mb{c}^{\top}\mb{x} \\
\mbox{subject to} && \mb{a}_i^0 \mb{x} \leq b_i - \rho ||\Delta_i \mb{x}||_2, \quad
\forall i=1,\dots,m.
\end{eqnarray*}

\end{theorem}
The intuition is as follows: for the case of ellipsoidal
uncertainty, the subproblem $\max_{\{\mb{a}_i \in \mathcal{U}_i\}}
\mb{a}_i^{\top} \mb{x} \leq b_i$, $\forall \, i$, is an
optimization over a quadratic constraint. The dual, therefore,
involves quadratic functions, which leads to
the resulting SOCP. \\ \\
{\bf Polyhedral Uncertainty:} Polyhedral uncertainty can be viewed
as a special case of ellipsoidal uncertainty
\cite{BenTalNemirovski99}. When $\mathcal{U}$ is
polyhedral, the subproblem becomes linear, and the robust
counterpart is equivalent to a linear optimization problem. To
illustrate this, consider the problem:
$$
\begin{array}{rl}
\min: & \mb{c}^{\top}\mb{x}  \\
\st: & \max_{\{\v{D}_i\v{a}_i\leq\v{d}_i\}}\mb{a}_i^{\top} \mb{x}
\leq b_i, \quad i=1,\dots,m.
\end{array}
$$
The dual of the subproblem (recall that $\mb{x}$ is not a variable
of optimization in the inner max) becomes:
$$
\left[ \begin{array}{rl}
\max: & \mb{a}_i^{\top} \mb{x} \\
\st: & \mb{D}_i \mb{a}_i \leq \mb{d}_i \end{array} \right]
\longleftrightarrow \left[ \begin{array}{rl}
\min: & \mb{p}_i^{\top} \mb{d}_i \\
\st: & \mb{p}_i^{\top}\mb{D}_i= \mb{x} \\
& \mb{p}_i \geq 0. \end{array} \right]
$$
and therefore the robust linear optimization now becomes:
$$
\begin{array}{rl}
\min: & \mb{c}^{\top} \mb{x} \\
\st: & \mb{p}_i^{\top} \mb{d}_i \leq b_i, \quad i=1,\dots,m \\
& \mb{p}_i^{\top}\mb{D}_i = \mb{x}, \quad i=1,\dots,m \\
& \mb{p}_i \geq 0, \quad i=1,\dots,m.
\end{array}
$$
Thus the size of such problems grows polynomially in the size of the
nominal problem and the dimensions of the uncertainty set. \\ \\
{\bf Cardinality Constrained Uncertainty:} Bertsimas and Sim
(\cite{BertsimasSim04a}) define a family of polyhedral uncertainty sets that encode a budget of uncertainty in terms of cardinality constraints:
the number of parameters of the problem that are allowed to vary from their nominal values. The uncertainty set $\mathcal{R}^D$ from our introductory example, is an instance of this. More generally, given an uncertain matrix, $\mb{A} = (a_{ij})$, suppose each component $a_{ij}$ lies in $[a_{ij} -
\hat{a}_{ij},a_{ij}+\hat{a}_{ij}]$. Rather than protect against the
case when every parameter can deviate, as in the original model of
Soyster (\cite{Soyster73}), we allow at most $\Gamma_i$ coefficients of row $i$
to deviate. Thus the positive number $\Gamma_i$
denotes the budget of uncertainty for the $i^{th}$
constraint, and just as with the ellipsoidal sizing, controls the trade-off
between the optimality of the solution, and its robustness to
parameter perturbation.\footnote{For the full details see
\cite{BertsimasSim04a}.} Given values $\Gamma_1,\dots,\Gamma_m$, the
robust formulation becomes: \be \label{eq:BSrobust}
\begin{array}{rll}
\min: & \mb{c}^{\top}\mb{x} \\
\st: & \sum_j a_{ij}x_j + \max_{\{S_i \subseteq J_i \,:\, |S_i| =
\Gamma_i\}}
\sum_{j \in S_i} \hat{a}_{ij}y_j \leq b_i & 1 \leq i \leq m \\
& -y_j \leq x_j \leq y_j & 1 \leq j \leq n \\
& \mb{l} \leq \mb{x} \leq \mb{u} & \\
& \mb{y} \geq \mb{0}. & \end{array}
\end{equation}
Because of the set-selection in the inner maximization, this problem is nonconvex. However, one can show that the natural convex relaxation is exact. Thus, relaxing and taking the dual of the inner maximization problem, one can show that the above is equivalent to the following linear formulation, and
therefore is tractable (and moreover is a linear optimization problem):
$$
\begin{array}{rll}
\max: & \mb{c}^{\top} \mb{x} & \\
\st: & \sum_j a_{ij} x_j + z_i\Gamma_i + \sum_j p_{ij} \leq b_i &
\forall \, i \\
& z_i + p_{ij} \geq \hat{a}_{ij} y_j & \forall \, i,j \\
& -y_j \leq x_j \leq y_j & \forall \, j \\
& \mb{l} \leq \mb{x} \leq \mb{u} & \\
& \mb{p} \geq \mb{0} & \\
& \mb{y} \geq \mb{0}. & \end{array}
$$
{\bf Norm Uncertainty:} Bertsimas et al.
\cite{BertsimasPachamanovaSim2004} show that robust linear
optimization problems with uncertainty sets described by more
general norms lead to convex problems with constraints related to
the dual norm. Here we use the notation $\mathrm{vec}(\v{A})$ to
denote the vector formed by concatenating all of the rows of the
matrix $\v{A}$.

\begin{theorem}\label{bsimcomp}
(Bertsimas et al., \cite{BertsimasPachamanovaSim2004}) With the
uncertainty set
\begin{eqnarray*}
\mathcal{U} &=& \{\v{A} \ | \ \|\v{M}(\mathrm{vec}(\v{A}) -
\mathrm{vec}(\v{\bar{A}}))\| \leq\Delta\},
\end{eqnarray*}
where $\v{M}$ is an invertible matrix, $\v{\bar{A}}$ is any constant
matrix, and $\|\cdot\|$ is any norm, Problem (\ref{roblpbg}) is
equivalent to the problem
\begin{eqnarray*}
\mathrm{minimize} && \v{c}^{\top}\v{x} \\
\mathrm{subject \ to} &&
\v{\bar{A}}_i^{\top}\v{x}+\Delta\|(\v{M}^{\top})^{-1}\v{x}_i\|^*\leq
b_i, \qquad i=1,\ldots,m,
\end{eqnarray*}
where $\v{x}_i\in\field{R}{(m\cdot n)\times 1}$ is a vector that
contains $\v{x}\in\field{R}{n}$ in entries $(i-1)\cdot n+1$ through
$i\cdot n$ and $0$ everywhere else, and $\|\cdot\|^*$ is the
corresponding dual norm of $\|\cdot\|$.
\end{theorem}

Thus the norm-based model shown in Theorem \ref{bsimcomp} yields an
equivalent problem with corresponding dual norm constraints. In particular,
the $\textit{l}_1$ and $\textit{l}_{\infty}$ norms result in linear
optimization problems, and the $\textit{l}_2$ norm results in a
second-order cone problem.

In short, for many choices of the uncertainty set, robust linear
optimization problems are tractable.

\subsection{\textbf{Robust quadratic optimization}}
\label{robqpsssec}
\emph{Quadratically constrained quadratic programs} (QCQP) have
defining functions $f_i(\v{x},\v{u}_i)$ of the form
\begin{eqnarray*}
f_i(\v{x},\v{u}_i) &=& \|\v{A}_i\v{x}\|^2+\v{b}_i^{\top}\v{x}+c_i.
\end{eqnarray*}
Second order cone programs (SOCPs) have
\begin{eqnarray*}
f_i(\v{x},\v{u}_i) &=& \|\v{A}_i\v{x}+\v{b}_i\|-\v{c}_i^{\top}\v{x}-d_i.
\end{eqnarray*}
For both classes, if the uncertainty set $\mathcal{U}$ is a single ellipsoid (called {\it simple ellipsoidal uncertainty}) the robust counterpart is a semidefinite optimization problem (SDP). If $\mathcal{U}$ is polyhedral or the intersection of ellipsoids, the robust counterpart is NP-hard (Ben-Tal and
Nemirovski, \cite{BenTalNemirovski98, BenTalNemirovski99,BenTalNemirovskiRoos02,BertsimasSim05}).

Following \cite{BenTalNemirovskiRoos02}, we illustrate here only how to obtain the explicit reformulation of
a robust quadratic constraint, subject to simple ellipsoidal
uncertainty. Consider the quadratic constraint
\begin{equation}
\label{eq:quadraticconstraint} \mb{x}^{\top}
\mb{A}^{\top}\mb{A}\mb{x} \leq 2\mb{b}^{\top} \mb{x} + \mb{c}, \quad
\forall (\mb{A},\mb{b},\mb{c}) \in \mathcal{U},
\end{equation}
where the uncertainty set $\mathcal{U}$ is an ellipsoid about a
nominal point $(\mb{A}^0,\mb{b}^0,\mb{c}^0)$:
$$
\mathcal{U} \bydef \left\{ (\mb{A},\mb{b},\mb{c}) :=
(\mb{A}^0,\mb{b}^0,\mb{c}^0) + \sum_{l=1}^L \mb{u}_l
(\mb{A}^l,\mb{b}^l,\mb{c}^l) \,:\, ||\mb{u}||_2 \leq 1 \right\}.
$$
As in the previous section, a vector $\mb{x}$ is feasible for the
robust constraint (\ref{eq:quadraticconstraint}) if and only if it
is feasible for the constraint:
$$
\left[ \begin{array}{rl} \max: & \mb{x}^{\top}
\mb{A}^{\top}\mb{A}\mb{x} - 2\mb{b}^{\top} \mb{x} - \mb{c} \\
\st: & (\mb{A},\mb{b},\mb{c}) \in \mathcal{U} \end{array} \right]
\leq 0.
$$
This is the maximization of a convex quadratic objective (when the
variable is the matrix $\mb{A}$, $\mb{x}^{\top}
\mb{A}^{\top}\mb{A}\mb{x}$ is quadratic and convex in $\mb{A}$ since
$\mb{x}\mb{x}^{\top}$ is always semidefinite) subject to a single
quadratic constraint. It is well-known that while this problem is
not convex (we are maximizing a convex quadratic) it nonetheless
enjoys a hidden convexity property (for an early reference, see
Brickman \cite{Brickman1961}) that allows it to be reformulated as a
(convex) semidefinite optimization problem. This is related to the so-called $S$-lemma (or $S$-procedure) in control
(e.g., Boyd et al. \cite{BoydElGhaouiFeronBalakrishnan1994}, P\'{o}lik and Terlaky \cite{PolikTerlaky07}):


The $S$-lemma essentially gives the boundary between
what we can solve exactly, and where solving the subproblem becomes
difficult. If the uncertainty set is an intersection of
ellipsoids or polyhedral, then exact solution of the subproblem is
NP-hard.\footnote{Nevertheless, there are some approximation results
available: \cite{BenTalNemirovskiRoos02}.}

Taking the dual of the SDP resulting from the $S$-lemma, we have an exact, convex reformulation of the subproblem in the RO
problem.

\begin{theorem} Given a vector $\mb{x}$, it is feasible to the
robust constraint (\ref{eq:quadraticconstraint}) if and only if
there exists a scalar $\tau \in \R$ such that the following matrix
inequality holds:
$$
\left(\begin{array}{c|ccc|c} c^0 + 2\mb{x}^{\top} \mb{b}^0 - \tau &
\frac{1}{2} c^1 + \mb{x}^{\top}
\mb{b}^1 &  \cdots & c^L + \mb{x}^{\top} \mb{b}^L & (\mb{A}^0 \mb{x})^{\top} \\
\hline
\frac{1}{2} c^1 + \mb{x}^{\top} \mb{b}^1 & \tau & & & (\mb{A}^1\mb{x})^{\top} \\
\vdots & & \ddots & & \vdots \\
\frac{1}{2} c^L + \mb{x}^{\top} \mb{b}^L & & & \tau & (\mb{A}^L\mb{x})^{\top} \\
\hline \mb{A}^0\mb{x} & \mb{A}^1\mb{x} & \cdots & \mb{A}^L\mb{x} & I
\end{array} \right) \succeq \mb{0}.
$$
\end{theorem}

\subsection{\textbf{Robust Semidefinite
Optimization}}\label{robsdpsssec}

With ellipsoidal uncertainty sets, robust counterparts of
semidefinite optimization problems are, in general, NP-hard (Ben-Tal and Nemirovski, \cite{BenTalNemirovski98}, Ben-Tal et al.
\cite{BenTalElGhaouiNemirovski00}). Similar negative results hold
even in the case of polyhedral uncertainty sets (Nemirovski,
\cite{Nemirovski93}). One exception (Boyd et al. \cite{BoydElGhaouiFeronBalakrishnan1994}) is when the uncertainty set is represented as \emph{unstructured norm-bounded uncertainty}. Such uncertainty takes the form
\begin{eqnarray*}
\v{A}_0(\v{x}) + \v{L}'(\v{x})\v{\zeta}\v{R}(\v{x}) + \v{R}(\v{x})\v{\zeta}\v{L}'(\v{x}),
\end{eqnarray*}
where $\v{\zeta}$ is a matrix with norm satisfying $\|\v{\zeta}\|_{2,2}\leq 1$, $\v{L}$ and $\v{R}$ are affine in the decision variables $\v{x}$, and at least one of $\v{L}$ or $\v{R}$ is independent of $\v{x}$.

In the general case, however, robust SDP is an intractable problem. Computing approximate solutions, i.e., solutions that are robust \emph{feasible} but not robust \emph{optimal} to robust semidefinite optimization problems has, as a consequence, received considerable attention (e.g., \cite{ElGhaouiOustryLebret98}, \cite{BenTalNemirovskiSDP01,
BenTalNemLMI02}, and \cite{BertsimasSim05}). These methods provide
bounds by developing inner approximations of the feasible set. The
goodness of the approximation is based on a measure of how close the inner approximation is to the true feasible set. Precisely, the measure for this is:
\begin{eqnarray*}
\rho(\text{AR}:\text{R}) &=& \inf\left\{\rho\geq 1 \ | \
X(\text{AR})\supseteq X(\mathcal{U}(\rho))\right\},
\end{eqnarray*}
where $X(\text{AR})$ is the feasible set of the approximate robust
problem and $X(\mathcal{U}(\rho))$ is the feasible set of the
original robust SDP with the uncertainty set ``inflated'' by a
factor of $\rho$. When the uncertainty set has ``structured norm bounded'' form, Ben-Tal and Nemirovski \cite{BenTalNemirovskiSDP01} develop an inner approximation such that
$\rho(\text{AR}:\text{R}) \leq \pi\sqrt{\mu}/2$, where $\mu$ is the
maximum rank of the matrices describing $\mathcal{U}$.

There has recently been additional work on Robust Semidefinite Optimization, for example exploiting sparsity \cite{OishiIsaka2009}, as well as in the area of control \cite{ElGhaouiNicolescu,CalafiorePolyak2001}.

\subsection{\textbf{Robust discrete
optimization}}\label{robintsssec}

Kouvelis and Yu \cite{KouvelisYu97} study robust models for some
discrete optimization problems, and show that the robust
counterparts to a number of polynomially solvable combinatorial
problems are NP-hard. For instance, the problem of minimizing the
maximum shortest path on a graph with only two scenarios for the
cost vector can be shown to be an NP-hard problem
\cite{KouvelisYu97}.

Bertsimas and Sim \cite{BertsimasSim03}, however, present a model
for cost uncertainty in which each coefficient $c_j$ is allowed to
vary within the interval $[\bar{c}_j,\bar{c}_j+d_j]$, with no more
than $\Gamma\geq 0$ coefficients allowed to vary. They then apply
this model to a number of combinatorial problems, i.e., they attempt
to solve
\begin{eqnarray*}
\text{minimize} && \bar{\v{c}}^{\top}\v{x}+ \max\limits_{\{S \ | \
S\subseteq N, \ |S|\leq\Gamma\}}\sum\limits_{j\in S}d_jx_j \\
\text{subject to} && \v{x}\in X,
\end{eqnarray*}
where $N=\{1,\ldots,n\}$ and $X$ is a fixed set. Under this model for uncertainty, the robust version of a
combinatorial problem may be solved by solving no more than $n+1$
instances of the underlying, nominal problem. This result extends to approximation algorithms for combinatorial
problems. For network flow problems, the above model
can be applied and the robust solution can be computed by solving a
logarithmic number of nominal, network flow problems.

Atamt\"{u}rk \cite{Atamturk05} shows that, under an appropriate
uncertainty model for the cost vector in a mixed 0-1 integer
program, there is a tight, linear programming formulation of the
robust problem with size polynomial in the size of a tight
linear programming formulation for the nominal problem.



\section{Choosing Uncertainty Sets}
\label{sec:probability}
In addition to tractability, a central question in the Robust Optimization literature has been probability guarantees on
feasibility under particular distributional assumptions for the
disturbance vectors. Specifically, what does robust feasibility
imply about probability of feasibility, i.e., what is the smallest
$\epsilon$ we can find such that
\begin{eqnarray*}
\v{x}\in X(\mathcal{U}) \Rightarrow \prob{f_i(\v{x},\v{u}_i)>0} \leq
\epsilon,
\end{eqnarray*}
under (ideally mild) assumptions on a distribution for $\v{u}_i$?

Such implications may be used as guidance for selection of a parameter representing the size of the uncertainty set. More generally, there are fundamental connections between distributional ambiguity, measures of risk, and uncertainty sets in robust optimization. In this section, we briefly discuss some of the connections in this vein.

\subsection{Probability Guarantees}

Probabilistic constraints, often called chance constraints in the literature, have a long history in stochastic optimization. Many approaches have been considered to address the computational challenges they pose (\cite{Prekopa95,NemirovskiShapiro05}), including work using sampling to approximate the chance constraints \cite{CalafioreCampi2003,CalafioreCampi06,ErdoganIyengar04}.

One of the early discussions of probability guarantees in RO traces back to Ben-Tal and Nemirovski \cite{BenTalNemirovski00}, who propose a robust model based on ellipsoids of radius $\Omega$ in the context of robust LP. Under this model, if the uncertain coefficients have bounded, symmetric support, they show that the corresponding robust feasible solutions must satisfy the constraint with high probability. Specifically, consider a linear constraint $\sum_j \tilde{a}_{ij}x_j\leq b_i$, where the coefficients $\tilde{a}_{ij}$ are uncertain and given by $\tilde{a}_{ij}=(1+\epsilon\xi_{ij})a_{ij}$, where $a_{ij}$ is a ``nominal'' value for the coefficient and $\{\xi_{ij}\}$ are zero mean, independent over $j$, and supported on $[-1,1]$. Then a robust constraint of the form
\begin{eqnarray*}
\sum\limits_j a_{ij}x_j + \epsilon\Omega\sqrt{\sum\limits_j a_{ij}^2x_j^2} \leq b_i^+,
\end{eqnarray*}
implies the robust solution satisfies the constraint with probability at least $1-e^{-\Omega^2/2}$. This bound holds for any such distribution on the finite support.

In a similar spirit, Bertsimas and Sim
\cite{BertsimasSim04a} propose an uncertainty set of the form
\begin{eqnarray}\label{ubsim}
\mathcal{U}_{\Gamma} &=& \left\{\v{\bar{A}} + \sum\limits_{j\in
J}z_j\hat{a}_j \ \Bigg| \ \|\v{z}\|_{\infty}\leq 1, \
\sum\limits_{j\in J}\mathbf{1}(z_j)\leq\Gamma\right\},
\end{eqnarray}
for the coefficients $\v{a}$ of an uncertain, linear constraint.
Here, $\mathbf{1}:\field{R}{}\rightarrow\field{R}{}$ denotes the
indicator function of $y$, i.e., $\mathbf{1}(y)=0$ if and only if
$y=0$, $\v{\bar{A}}$ is a vector of ``nominal'' values,
$J\subseteq\{1,\ldots,n\}$ is an index set of uncertain
coefficients, and $\Gamma\leq |J|$ is an integer\footnote{The
authors also consider $\Gamma$ non-integer, but we omit this
straightforward extension for notational convenience.} reflecting
the number of coefficients which are allowed to deviate from their
nominal values. The dual formulation of this as a linear
optimization is discussed in Section \ref{sec:structuretractability}. The following then holds.
\begin{theorem}\label{bsimgamthm}
(Bertsimas and Sim \cite{BertsimasSim04a}) Let $\v{x}^*$ satisfy the constraint
\begin{eqnarray*}
\max\limits_{\v{a}\in\mathcal{U}_{\Gamma}}\v{a}^{\top}\v{x}^*\leq b,
\end{eqnarray*}
where $\mathcal{U}_{\Gamma}$ is as in (\ref{ubsim}). If the random
vector $\tilde{\v{a}}$ has independent components with $a_j$
distributed symmetrically on
$[\bar{a}_j-\hat{a}_j,\bar{a}_j+\hat{a}_j]$ if $j\in J$ and
$a_j=\bar{a}_j$ otherwise, then
\begin{eqnarray*}
\prob{\tilde{\v{a}}^{\top}\v{x}^*>b}\leq e^{-\frac{\Gamma^2}{2|J|}}.
\end{eqnarray*}
\end{theorem}

In the case of linear optimization with only partial moment
information (specifically, known mean and covariance), Bertsimas et
al. \cite{BertsimasPachamanovaSim2004} prove guarantees for the
general norm uncertainty model used in Theorem \ref{bsimcomp}. For
instance, when $\|\cdot\|$ is the Euclidean norm, and $\v{x}^*$ is feasible to the robust problem, Theorem \ref{bsimcomp} can be shown
\cite{BertsimasPachamanovaSim2004} to imply the guarantee
\begin{eqnarray*}
\prob{\tilde{\v{a}}^{\top}\v{x}^*>b}\leq \frac{1}{1+\Delta^2},
\end{eqnarray*}
where $\Delta$ is the radius of the uncertainty set, and the mean
and covariance are used for $\v{\bar{A}}$ and $\v{M}$, respectively.

For more general robust conic optimization problems, results on
probability guarantees are more elusive. Bertsimas and Sim are able
to prove probability guarantees for their approximate robust
solutions in \cite{BertsimasSim05}. In Chen et al. \cite{chensimsun06}, more general deviation measures are considered that capture distributional skewness, leading to improved probability guarantees. Also of interest is the work of Paschalidis and Kang on probability guarantees and uncertainty set selection when the entire distribution is available \cite{PaschalidisKang2005}.

\subsection{Distributional Uncertainty}
The issue of limited distributional information is central and has been the subject of considerable research in the decision theory literature. This work closely connects to robustness considerations and provides potential guidance and economic meaning to the choice of particular uncertainty sets.

Consider a function $u(\v{x},\xi)$ where $\xi$ is a random parameter on some measure space $(\Omega,\mathcal{F})$. For the purposes of this discussion, let $u$ be a concave, nondecreasing payoff function. In many situations, it may be unreasonable to expect the decision maker to have a full description of the distribution of $\xi$, but instead knows the distribution to be confined to some set of distributions $\mathcal{Q}$. Using a well-known duality result that traces back to at least the robust statistics literature (e.g., Huber \cite{Huber81}), one can establish that for any set $\mathcal{Q}$, there exists a convex, non-increasing, translation-invariant, positive homogeneous function $\mu$ on the induced space of random variables, such that
\begin{eqnarray}\label{distrob}
\inf\limits_{\mathbb{Q}\in\mathcal{Q}} \ev{\mathbb{Q}}{u(\v{x},\xi)} \ \geq \ 0 &\Leftrightarrow& \mu(u(\v{x},\xi)) \ \leq \ 0.
\end{eqnarray}
The function in this representation falls precisely into the class of \emph{coherent risk measures} popularized by Artzner et al. \cite{ArtznerDelbaenEberHeath99}. These functions provide an economic interpretation in terms of a capital requirement: if $X$ is a random variable (e.g., return), $\mu(X)$ represents the amount of money required to be added to $X$ in order to make it ``acceptable,'' given utility function $u$. The properties listed above are natural in a risk management setting: monotonicity states that one position that always pays off more than another should be deemed less risky; translation invariance means the addition of a sure amount to a position reduces the risk by precisely that amount; positive homogeneity means risks scale equally with the size of the stakes; and convexity means diversification among risky positions should be encouraged.

The above observation implies an immediate connection between these risk management tools, distributional ambiguity, and robust optimization. These connections have been explored in recent work on robust optimization. Natarajan et al. \cite{NatarajanPachamanoveSim05} investigate this connection with a focus on inferring risk measures from uncertainty sets.

Bertsimas and Brown \cite{BertsimasBrown05a} examine the question from the opposite perspective: namely, with risk preferences specified by a coherent risk measure, they examine the implications for uncertainty set structure in robust linear optimization problems. Due to the duality above, a risk constraint of the form $\mu(\tilde{\v{a}}'\v{x}-b)\leq 0$ on a linear constraint with an uncertain vector $\tilde{\v{a}}$ can be equivalently expressed as
\begin{eqnarray*}
\v{a}'\v{x} &\geq& b \qquad \forall \ \v{a}\in\mathcal{U},
\end{eqnarray*}
where $\mathcal{U} = \ch{\{\ev{\mathbb{Q}}{\v{a}} \ : \ \mathbb{Q}\in\mathcal{Q}\}}$
and $\mathcal{Q}$ is the generating family for $\mu$.

For a concrete application of this, one of most famous coherent risk measures is the \emph{conditional value-at-risk} (CVaR), defined as
\begin{eqnarray*}
\mu(X) &\triangleq& \inf\limits_{\nu\in\field{R}{}}\left\{ \nu +
\frac{1}{\alpha} \ev{}{(-\nu-X)^+}\right\},
\end{eqnarray*}
for any $\alpha\in (0,1]$. For atomless distributions, CVaR is equivalent to the expected value of the random variable conditional on it being in its lower $\alpha$ quantile.

Consider the case when the uncertain vector $\tilde{\v{a}}$ follows a discrete distribution with support $\{\v{a}_1,\ldots,\v{a}_N\}$ and corresponding probabilities $\{\v{p}_1,\ldots,\v{p}_N\}$. The generating family for CVaR in this case is $\mathcal{Q}=\{\v{q}\in\Delta^N \ : \ q_i\leq p_i/\alpha\}$. This leads to the uncertainty set
\begin{eqnarray*}
\mathcal{U} &=& \ch{ \left\{\frac{1}{\alpha}\sum\limits_{i\in
I}p_i\v{a}_i+ \left( 1-\frac{1}{\alpha}\sum\limits_{i\in
I}p_i\right) \v{a}_j \ : \ I\subseteq \{1,\ldots,N\}, \
j\in\{1,\ldots,N\}\setminus I, \ \sum\limits_{i\in
I}p_i\leq\alpha\right\}}.
\end{eqnarray*}
This set is a polytope, and therefore the robust optimization problem in this case may be reformulated as a linear program. When $p_i=1/N$ and $\alpha=j/N$ for some $j\in\field{Z}{}_+$, this has the interpretation of the convex hull of all $j$-point averages of $\mathcal{A}$.

Despite its popularity, CVaR represents only a special case of a much broader class of coherent risk measures that are \emph{comonotone}. These risk measures satisfy the additional property that risky positions that ``move together'' in all states cannot be used to hedge one another. Extending a result from Dellacherie \cite{Dellacherie}, Schmeidler \cite{Schmeidler86} shows that the class of such risk measures is precisely the same as the set of functions representable as \emph{Choquet integrals} (Choquet, \cite{Choquet54}). Choquet integrals are the expectation under a set function that is non-additive and are a classical approach towards dealing with ambiguous distributions in decision theory. Bertsimas and Brown \cite{BertsimasBrown05b} discuss how one can form uncertainty sets in RO with these types of risk measures on discrete event spaces.

The use of a discrete probability space may be justified in situations when samples of the uncertainty are available. Delage and Ye \cite{delageye09} have proposed an approach to the distribution-robust problem
\begin{eqnarray*}
\text{minimize}_{\v{x}\in X} \max\limits_{f_{\xi}\in\mathcal{D}} \ev{\xi}{h(\v{x},\xi)},
\end{eqnarray*}
where $\xi$ is a random parameter with distribution $f_{\xi}$ on some set of distributions $\mathcal{D}$ supported on a bounded set $\mathcal{S}$, $h$ is convex in the decision variable $\v{x}$, and $X$ is a convex set. They consider sets of distributions $\mathcal{D}$ based on moment uncertainty with a particular focus on sets that have uncertainty in the mean and covariance of $\xi$. They then consider the problem when one has independent samples $\xi_1,\ldots,\xi_M$ and focus largely on the set
\begin{eqnarray*}
&& \mathcal{D}_1(\mathcal{S},\hat{\v{\mu}}_0,\hat{\v{\Sigma}}_0,\gamma_1,\gamma_2) \triangleq \\
 && \left\{ \prob{\xi\in S} = 1 \ : \ (\ev{}{\xi}-\hat{\v{\mu}}_0)'\hat{\v{\Sigma}_0}^{-1}(\ev{}{\xi}-\hat{\v{\mu}}_0)\leq \gamma_1, \ \ev{}{(\xi-\hat{\v{\mu}}_0)'(\xi-\hat{\v{\mu}}_0)}\preceq \gamma_2\hat{\v{\Sigma}}_0\right\}.
\end{eqnarray*}
The above problem can be solved in polynomial time, and, with proper choices of $\gamma_1,\gamma_2$ and $M$, the resulting optimal value provides an upper bound on the expected cost with high probability. In the case of $h$ as a piecewise linear, convex function, the resulting problem reduces to solving an SDP. This type of approach seems highly practical in settings (prevalent in many applications, e.g., finance) where samples are the only relevant information a decision maker has on the underlying distribution.

Related to distributional uncertainty is the work in \cite{XuCaramanisMannorDR}. Here, Xu, Caramanis and Mannor show that any robust optimization problem is equivalent to a distributionally robust problem. Using this equivalence to robust optimization, they show how robustness can guarantee consistency in sampled problems, even when the nominal sampled problem fails to be consistent.

More general types of robust optimization models have been explored, and such approaches draw further connections to research in decision theory. Ben-Tal et al. \cite{BentalBertsimasBrown08} propose an approach called \emph{soft robust optimization} applicable in settings of distributional ambiguity. They modify the constraint (\ref{distrob}) and consider the more general constraint
\begin{eqnarray*}
\inf\limits_{\mathbb{Q}\in\mathcal{Q}(\epsilon)} \ev{\mathbb{Q}}{f(\v{x},\xi)} & \geq& -\epsilon \qquad\forall\epsilon\geq 0,
\end{eqnarray*}
where $\{\mathcal{Q}(\epsilon)\}_{\epsilon\geq 0}$ is a set of sets of distributions, nondecreasing and convex on $\epsilon\geq 0$. This set of constraints considers different sized uncertainty sets with increasingly looser feasibility requirements as the uncertainty size grows; as such, it provides a potentially less conservative approach to RO than (\ref{distrob}). This approach connects to the approach of \emph{convex risk measures} (F\"{o}llmer and Schied, \cite{FollmerSchied02}), a generalization of the coherent risk measures mentioned above. Under a particular form for $\mathcal{Q}(\epsilon)$ based on relative entropy deviations, this model recovers the \emph{multiplier preferences} of Hansen and Sargent \cite{hansensargent01}, who develop their approach from robust control ideas in order to deal with model mis-specification in the decision making of economic agents (see also Maccheroni et al. \cite{maccheroni06} for a generalization known as \emph{variational preferences}).

In short, there has been considerable work done in the domain of uncertainty set construction for RO. Some of this work focuses on the highly practical matter of implied probability guarantees under mild distributional assumptions or under a sufficiently large number of samples; other work draws connections to objects that have been axiomatized and developed in the decision theory literature over the past several decades.

\section{Robust Adaptable Optimization}\label{adaptsec}

Thus far this paper has addressed optimization in the static, or
one-shot case: the decision-maker considers a single-stage
optimization problem affected by uncertainty. In this formulation,
all the decisions are implemented simultaneously, and in particular,
before any of the uncertainty is realized. In dynamic (or sequential)
decision-making problems this single-shot assumption is restrictive and
conservative. For example, in the inventory control example we discuss below, this would
correspond to making all ordering decisions up front, without flexibility to adapt
to changing demand patterns.

Sequential decision-making appears in a broad range of applications in many areas of
engineering and beyond. There has been extensive work in optimal and robust control (e.g., the textbooks
 \cite{DullerudPaganini1999,ZhouDoyleGlover1996}, or the articles \cite{FanTitsDoyle1991,GoulartKerrigan2005,GriederParriloMorari2003,KerriganMaciejowski2003}, and references therein), and approximate and exact dynamic programming (e.g., see the textbooks
  \cite{Bertsekas95dp1,BertsekasTsitsiklis1996,Puterman1994}). In this section, we consider
modeling approaches to incorporate sequential decision-making into the robust optimization
framework.

\subsection{Motivation and Background}
In what follows, we refer to the {\it static} solution as the case
where the $\mb{x}_i$ are all chosen at time 1 before any
realizations of the uncertainty are revealed. The {\it dynamic}
solution is the fully adaptable one, where $\mb{x}_i$ may have
arbitrary functional dependence on past realizations of the
uncertainty.

The question as to when adaptability has value is an interesting one that has received some attention. The papers by Dean, Goemans and Vondr\'{a}k (\cite{DeanGoemansVondrak2004,GoemansVondrak2006}) consider the value of adaptability in the context of stochastic optimization problems. They show there that for the stochastic knapsack problem, the value of adaptability is bounded: the value of the optimal adaptive solution is no more than a constant factor times the value of the optimal non-adaptive solution. In \cite{BertsimasGoyal2009}, Bertsimas and Goyal consider a two-stage mixed integer stochastic optimization problem with uncertainty in the right-hand-side. They show that a static robust solution approximates the fully-adaptable two-stage solution for the stochastic problem to within a factor of two, as long as the uncertainty set and the underlying measure, are both symmetric.

Despite the results for these cases, we would generally expect approaches that explicitly incorporate adaptivity to substantially outperform static approaches in multi-period problems. There are a number of approaches.


\subsubsection*{Receding Horizon}
The most straightforward extension of the single-shot Robust
Optimization formulation to that of sequential decision-making, is
the so-called receding horizon approach. In this formulation, the
static solution over all stages is computed, and the first-stage
decision is implemented. At the next stage, the process is repeated. In the control literature this is known as open-loop feedback. While this approach is typically tractable, in many cases it may be far
from optimal. In particular, because it is computed without any
adaptability, the first stage decision may be overly conservative.

\subsubsection*{Stochastic Optimization}
In Stochastic Optimization, the basic problem of interest is the
so-called complete recourse problem (for the
basic definitions and setup, see \cite{BirgeLouveaux97,
Infanger1994, Prekopa95}, and references therein). In this setup,
the feasibility constraints of a single-stage Stochastic
Optimization problem are relaxed and moved into the objective
function by assuming that after the first-stage decisions are
implemented and the uncertainty realized, the decision-maker has
some recourse to ensure that the constraints are satisfied. The
canonical example is in inventory control where in case of shortfall
the decision-maker can buy inventory at a higher cost (possibly from
a competitor) to meet demand. Then the problem becomes one of
minimizing expected cost of the two-stage problem. If there is no
complete recourse (i.e., not every first-stage decision can be
completed to a feasible solution via second-stage actions) and
furthermore the impact and cost of the second-stage actions are
uncertain at the first stage, the problem becomes considerably more
difficult. The feasibility constraint in particular is much more
difficult to treat, since it cannot be entirely brought into the
objective function.

When the uncertainty is assumed to take values in a finite set of
small cardinality, the two-stage problem is tractable, and even for
larger cardinality (but still finite) uncertainty sets (called
scenarios), large-scale linear programming techniques such as
Bender's decomposition can be employed to obtain a tractable
formulation (see, e.g., \cite{BertsimasTsitsiklis1997}). In the case
of incomplete recourse where feasibility is not guaranteed,
robustness of the first-stage decision may require a very large
number of scenarios in order to capture enough of the structure of
the uncertainty. In the next section, we discuss a robust adaptable
approach called Finite Adaptability that seeks to circumvent this
issue.

Finally, even for small cardinality sets, the multi-stage complexity
explodes in the number of stages (\cite{Shapiro05}). This is a
central problem of multi-stage optimization, in both the robust and
the stochastic formulations.

\subsubsection*{Dynamic Programming}
Sequential decision-making under uncertainty has traditionally fallen under the purview
of Dynamic Programming, where many exact and approximate techniques have been
developed -- we do not review this work here, but rather refer the reader to the
books \cite{Bertsekas95dp1},  \cite{BertsekasTsitsiklis1996}, and \cite{Puterman1994}.
The Dynamic Programming framework has
been extended to the robust Dynamic Programming and robust
MDP setting, where the payoffs and the dynamics are not exactly
known, in Iyengar \cite{Iyengar2005} and Nilim and El Ghaoui
\cite{NilimElghaoui2004a}, and then also in Xu and Mannor
\cite{XuMannor2006}. Dynamic Programming yields tractable
algorithms precisely when the Dynamic Programming recursion does not
suffer from the curse of dimensionality. As the papers cited above
make clear, this is a fragile property of any problem, and is
particularly sensitive to the structure of the uncertainty. Indeed,
the work in \cite{Iyengar2005, NilimElghaoui2004a, XuMannor2006,
DelageMannor2006} assumes a special property of the uncertainty set
(``rectangularity'') that effectively means that the decision-maker
gains nothing by having future stage actions depend explicitly on
past realizations of the uncertainty.

This section is devoted precisely to this problem: the dependence of future actions on past realizations of the uncertainty.

\subsection{Tractability of Robust Adaptable Optimization}
The uncertain multi-stage problem with deterministic set-based
uncertainty, i.e., the robust multi-stage formulation, was first
considered in \cite{BenTalGoryashkoGuslitzerNemirovski03}. There,
the authors show that the two-stage linear problem with
deterministic uncertainty is in general $NP$-hard. Consider the
generic two-stage problem: \be \label{eq:generic2stage}
\begin{array}{rl}
\min: & \mb{c}^{\top}\mb{x}_1 \\
\st: & \mb{A}_1(\mb{u}) \mb{x}_1 + \mb{A}_2(\mb{u})\mb{x}_2(\mb{u})
\leq \mb{b}, \quad \forall \mb{u} \in {\cal U}.
\end{array}
\ee Here, $\mb{x}_2(\cdot)$ is an arbitrary function of $\mb{u}$. We
can rewrite this explicitly in terms of the feasible set for the
first stage decision:
\begin{equation}
\label{eq:quantformulation}
\begin{array}{rl}
\min: & \mb{c}^{\top}\mb{x}_1 \\
\st: & \mb{x}_1 \in \left\{ \mb{x}_1 \,:\, \forall \mb{u} \in {\cal
U}, \, \exists \mb{x}_2 \,\, \st \mb{A}_1(\mb{u})\mb{x}_1 +
\mb{A}_2(\mb{u}) \mb{x}_2 \leq \mb{b} \right\}.
\end{array}
\end{equation}
The feasible set is convex, but nevertheless the optimization
problem is in general intractable. Consider a simple example given
in \cite{BenTalGoryashkoGuslitzerNemirovski03}: \be
\label{eq:NPhard2stage}
\begin{array}{rl}
\min: & x_1 \\
\st: & x_1 - \mb{u}^{\top} \mb{x}_2(\mb{u}) \geq 0 \\
& \mb{x}_2(\mb{u}) \geq  \mb{B}\mb{u} \\
& \mb{x}_2(\mb{u}) \leq  \mb{B}\mb{u}.
\end{array}
\ee It is not hard to see that the feasible first-stage decisions
are given by the set:
$$
\{x_1 \,:\, x_1 \geq \mb{u}^{\top}\mb{B}\mb{u},\quad \forall \mb{u}
\in {\cal U}\}.
$$
The set is, therefore, a ray in $\R^1$, but determining the left
endpoint of this line requires computing a maximization of a
(possibly indefinite) quadratic $\mb{u}^{\top}\mb{B}\mb{u}$, over
the set ${\cal U}$. In general, this problem is NP-hard
(see, e.g., \cite{GareyJohnson79}).

\subsection{Theoretical Results}
Despite the hardness result above, there has been
considerable effort devoted to obtaining different approximations
and approaches to the multi-stage optimization problem.

\subsubsection{Affine Adaptability}
\label{ssec:affine} In \cite{BenTalGoryashkoGuslitzerNemirovski03},
the authors formulate an approximation to the general robust
multi-stage optimization problem, which they call the {\it Affinely
Adjustable Robust Counterpart} (AARC). Here, they explicitly
parameterize the future stage decisions as affine functions of the
revealed uncertainty. For the two-stage problem
(\ref{eq:generic2stage}), the second stage variable,
$\mb{x}_2(\mb{u})$, is parameterized as:
$$
\mb{x}_2(\mb{u}) = \mb{Q}\mb{u} + \mb{q}.
$$
Now, the problem becomes:
$$
\begin{array}{rl}
\min: & \mb{c}^{\top}\mb{x}_1 \\
\st: & \mb{A}_1(\mb{u})\mb{x}_1 + \mb{A}_2(\mb{u})[ \mb{Q}\mb{u} +
\mb{q}] \leq \mb{b}, \quad \forall \mb{u} \in {\cal U}.
\end{array}
$$
This is a single-stage RO. The
decision-variables are $(\mb{x}_1,\mb{Q},\mb{q})$, and they are all
to be decided before the uncertain parameter, $\mb{u} \in {\cal U}$,
is realized.

In the generic formulation of the two-stage problem
(\ref{eq:generic2stage}), the functional dependence of
$\mb{x}_2(\cdot)$ on $\mb{u}$ is arbitrary. In the affine  formulation,
the resulting problem is a linear optimization
problem with uncertainty. The parameters of the problem, however,
now have a quadratic dependence on the uncertain parameter $\mb{u}$.
Thus in general, the resulting robust linear optimization will not
be tractable -- consider again the example
$(\ref{eq:NPhard2stage})$.

Despite this negative result, there are some positive complexity
results concerning the affine model. In order to present these, let
us explicitly denote the dependence of the optimization parameters,
$\mb{A}_1$ and $\mb{A}_2$, as:
$$
[\mb{A}_1,\mb{A}_2](\mb{u}) =
[\mb{A}_1^{(0)},\mb{A}_2^{(0)}]+\sum_{l=1}^{L} u_l
[\mb{A}_1^{(l)},\mb{A}_2^{(l)}].
$$
When we have $\mb{A}_2^{(l)} = \mb{0}$, for all $l \geq 1$, the
matrix multiplying the second stage variables is constant. This
setting is known as the case of {\it fixed recourse}.
We can now write the second stage variables explicitly in terms of the
columns of the matrix $\mb{Q}$. Letting $\mb{q}^{(l)}$ denote the
$l^{th}$ column of $\mb{Q}$, and $\mb{q}^{(0)} = \mb{q}$ the
constant vector, we have:
\begin{eqnarray*}
\mb{x}_2 &=& \mb{Q}\mb{u} + \mb{q}_0 \\
&=& \mb{q}^{(0)} + \sum_{l=1}^{L} u_l \mb{q}^{(l)}.
\end{eqnarray*}
Letting $\mb{\chi} = (\mb{x}_1,\mb{q}^{(0)},\mb{q}^{(1)}, \dots, \mb{q}^{(L)})$ denote
the full decision vector, we can write the $i^{th}$ constraint as
\begin{eqnarray*}
0 &\leq& (\mb{A}_1^{(0)}\mb{x}_1 + \mb{A}_2^{(0)}\mb{q}^{(0)} -
\mb{b})_i+ \sum_{l=1}^L u_l(\mb{A}_1^{(l)}\mb{x}_1 + \mb{A}_2
\mb{q}^{(l)})_i \\
&=& \sum_{l=0}^L a_l^i(\mb{\chi}),
\end{eqnarray*}
where we have defined
$$
a_l^i \bydef a_l^i(\mb{\chi}) \bydef (\mb{A}_1^{(l)}\mb{x}_1 +
\mb{A}_2^{(l)} \mb{q}^{(l)})_i, \qquad a_0^i \bydef
(\mb{A}_1^{(0)}\mb{x}_1 + \mb{A}_2^{(0)}\mb{q}^{(0)} - \mb{b})_i.
$$
\begin{theorem}[\cite{BenTalGoryashkoGuslitzerNemirovski03}]
Assume we have a two-stage linear optimization with fixed recourse,
and with conic uncertainty set:
$$
{\cal U} = \{\mb{u}\,:\, \exists \xi \, \st \, \mb{V}_1 \mb{u} +
\mb{V}_2 \mb{\xi} \geq_{{\cal K}} \mb{d}\} \subseteq \R^{L},
$$
where ${\cal K}$ is a convex cone with dual ${\cal K}^{\ast}$. If ${\cal U}$
has nonempty interior, then the AARC can be
reformulated as the following optimization problem:
\begin{eqnarray*}
\min: && \mb{c}^{\top} \mb{x}_1 \\
\st: && \mb{V}_1 \lambda^i -
a^i(\mb{x}_1,\mb{q}^{(0)},\dots,\mb{q}^{(L)}) = 0, \quad i=1,\dots,m
\\
&& \mb{V}_2 \lambda^i = 0, \quad i=1,\dots,m \\
&& \mb{d}^{\top} \lambda^i +
a_0^i(\mb{x}_1,\mb{q}^{(0)},\dots,\mb{q}^{(L)}) \geq 0, \quad
i=1,\dots,m \\
&& \lambda^i \geq_{{\cal K}^{\ast}} 0, \quad i=1,\dots,m.
\end{eqnarray*}
\end{theorem}
If the cone ${\cal K}$ is the positive orthant, then the AARC given
above is an LP. \\ \\
The case of non-fixed recourse is more difficult because of the
quadratic dependence on $\mb{u}$. Note that the example in
(\ref{eq:NPhard2stage}) above involves an uncertainty-affected
recourse matrix. In this non-fixed recourse case, the robust
constraints have a component that is quadratic in the uncertain
parameters, $\mb{u}_i$. These robust constraints then become:
$$
\left[ \mb{A}_1^{(0)} + \sum \mb{u}_l A_1^{(1)} \right] \mb{x}_1 +
\left[ \mb{A}_2^{(0)} + \sum \mb{u}_l A_2^{(1)} \right] \left[
\mb{q}^{(0)} + \sum \mb{u}_l \mb{q}^{(l)} \right] - \mb{b} \leq
\mb{0}, \quad \forall \mb{u} \in {\cal U},
$$
which can be rewritten to emphasize the quadratic dependence on
$\mb{u}$, as
$$
\left[ \mb{A}_1^{(0)} \mb{x}_1 + \mb{A}_2^{(0)} \mb{q}^{(0)} -
\mb{b} \right] + \sum \mb{u}_l \left[\mb{A}_1^{(l)}\mb{x}_1 +
\mb{A}_2^{(0)} \mb{q}^{(l)} + \mb{A}_2^{(l)} \mb{q}^{(0)} \right] +
\left[ \sum \mb{u}_k \mb{u}_l \mb{A}_2^{(k)} \mb{q}^{(l)} \right]
\leq 0, \quad \forall \mb{u} \in {\cal U}.
$$
Writing
\begin{eqnarray*}
\mb{\chi} &\bydef& (\mb{x}_1,\mb{q}^{(0)},\dots,\mb{q}^{(L)}), \\
\alpha_i(\mb{\chi}) &\bydef& -[\mb{A}_1^{(0)} \mb{x}_1 +
\mb{A}_2^{(0)}
\mb{q}^{(0)} - \mb{b}]_i \\
\beta_i^{(l)}(\mb{\chi}) &\bydef& -\frac{[\mb{A}_1^{(l)}\mb{x}_1 +
\mb{A}_2^{(0)} \mb{q}^{(l)} - \mb{b}]_i}{2}, \quad l=1,\dots,L \\
\Gamma_i^{(l,k)}(\mb{\chi}) &\bydef& - \frac{[\mb{A}_2^{(k)}
\mb{q}^{(l)} + \mb{A}_2^{(l)}\mb{q}^{(k)}]_i}{2}, \quad
l,k=1,\dots,L,
\end{eqnarray*}
the robust constraints can now be expressed as:
\begin{equation}
\label{eq:AARCconstraint}
 \alpha_i(\mb{\chi}) +
2\mb{u}^{\top} \beta_i(\mb{\chi}) + \mb{u}^{\top}
\Gamma_i(\mb{\chi}) \mb{u} \geq 0, \quad \forall \mb{u} \in {\cal
U}.
\end{equation}
\begin{theorem}[\cite{BenTalGoryashkoGuslitzerNemirovski03}]
Let our uncertainty set be given as the intersection of ellipsoids:
$$
{\cal U} \bydef \{\mb{u} \,:\, \mb{u}^{\top} (\rho^{-2}S_k) \mb{u}
\leq 1, \,\, k=1,\dots,K\},
$$
where $\rho$ controls the size of the ellipsoids.
Then the original AARC problem can be approximated by the following
semidefinite optimization problem:
\begin{equation}
\label{eq:AARCSDP}
\begin{array}{rl}
\min: & \mb{c}^{\top} \mb{x}_1 \\
\st: & \left(
\begin{array}{c|c} \Gamma_i(\mb{\chi}) + \rho^{-2}\sum_{k=1}^K \lambda_k
S_k & \beta_i(\mb{\chi}) \\ \hline \beta_i(\mb{\chi})^{\top} &
\alpha_i(\mb{\chi}) - \sum_{k=1}^K \lambda_k^{(i)}
\end{array} \right) \succeq \mb{0}, \,\, i=1,\dots, m \\
& \lambda^{(i)} \geq 0, \,\, i=1,\dots,m
\end{array}
\end{equation}
\end{theorem}

The constant $\rho$ in the definition of the uncertainty set ${\cal
U}$ can be regarded as a measure of the level of the uncertainty.
This allows us to give a bound on the tightness of the
approximation. Define the constant
$$
\gamma \bydef \sqrt{ 2\ln \left( 6 \sum_{k=1}^K \mbox{Rank}(S_k)
\right) } .
$$
Then we have the following.
\begin{theorem}[\cite{BenTalGoryashkoGuslitzerNemirovski03}] Let
${\cal X}_{\rho}$ denote the feasible set of the AARC with noise
level $\rho$. Let ${\cal X}_{\rho}^{\mbox{{\tiny approx}}}$ denote
the feasible set of the SDP approximation to the AARC with
uncertainty parameter $\rho$. Then, for $\gamma$ defined as above,
we have the containment:
$$
{\cal X}_{\gamma \rho} \subseteq {\cal X}_{\rho}^{\mbox{{\tiny
approx}}} \subseteq {\cal X}_{\rho}.
$$
\end{theorem}
\noindent This tightness result has been improved; see
\cite{DerinkuyuPinar2006}.

There have been a number of applications building upon affine
adaptability, in a wide array of areas:

\begin{enumerate}
\item Integrated circuit design: In \cite{ManiSinghOrshansky2006}, the affine
adjustable approach is used to model the yield-loss optimization in
chip design, where the first stage decisions are the pre-silicon
design decisions, while the second-stage decisions represent post-silicon
tuning, made after the manufacturing variability is
realized and can then be measured.
\item Comprehensive Robust Optimization: In
\cite{BenTalBoydNemirovski2006}, the authors extend the robust
static, as well as the affine adaptability framework, to soften the
hard constraints of the optimization, and hence to reduce the
conservativeness of robustness. At the same time, this controls the
infeasibility of the solution even when the uncertainty is realized
outside a nominal compact set. This has many applications, including
portfolio management, and optimal control.
\item Network flows and Traffic Management: In \cite{OrdonezZhao2005},
the authors consider the robust capacity expansion of a network flow
problem that faces uncertainty in the demand, and also the travel
time along the links. They use the adjustable framework of
\cite{BenTalGoryashkoGuslitzerNemirovski03}, and they show that for
the structure of uncertainty sets they consider, the resulting
problem is tractable. In \cite{MudchanatongsukOrdonezLiu2005}, the
authors consider a similar problem under transportation cost and
demand uncertainty, extending the work in \cite{OrdonezZhao2005}.
\item Chance constraints: In \cite{ChenSimSunZhang07}, the authors
apply a modified model of affine adaptability to the stochastic
programming setting, and show how this can improve approximations of chance constraints. In \cite{ErdoganIyengar05}, the
authors formulate and propose an algorithm for the problem of
two-stage convex chance constraints when the underlying distribution
has some uncertainty (i.e., an \emph{ambiguous} distribution).
\item Numerous other applications have been considered, including portfolio management \cite{CalafioreCampi2005b,TakedaTaguchiTutuncu04}, coordination in wireless networks \cite{YunCaramanis2008}, robust control \cite{GoulartKerriganMaciejowski06}, and model adaptive control.

\end{enumerate}
Additional work in affine adaptability has been done in \cite{ChenSimSunZhang07}, where the authors consider modified linear decision rules in the context of only partial distibutional knowledge, and within that framework derive tractable approximations to the resulting robust problems. See also Ben-Tal et al. \cite{BentalElGhaouiNemirovskiBook} for a detailed discussion of affine decision rules in multistage optimization. Recently, \cite{BertsimasIancuParrilo2009a} have given conditions under which affine policies are in fact optimal, and affine policies have been extended to higher order polynomial adaptability in \cite{BertsimasCaramanis2007, BertsimasIancuParrilo2009b}.

\subsubsection{Finite Adaptability}
The framework of Finite Adaptability, introduced in Bertsimas and
Caramanis \cite{BertsimasCaramanisTAC} and Caramanis
\cite{CaramanisPhD}, treats the discrete setting by modeling the second-stage variables, $\mb{x}_2(\mb{u})$, as piecewise
constant functions of the uncertainty, with $k$ pieces. One advantage of such an approach is that, due to the
inherent finiteness of the framework, the resulting formulation can
accommodate discrete variables. In addition, the level of
adaptability can be adjusted by changing the number of pieces in the piecewise constant second stage variables. (For an example from
circuit design where such second stage limited adaptability
constraints are physically motivated by design considerations, see
\cite{ManiCaramanisOrshansky2007,SinghHeCaramanisOrshansky2009}).

If the partition of the uncertainty set 
is fixed, then the resulting problem retains
the structure of the original nominal problem, and the number of
second stage variables grows by a factor of $k$.
In general, computing the optimal partition into even two regions
is NP-hard \cite{BertsimasCaramanisTAC}, however, if any one of the three
quantities: (a) dimension of the uncertainty; (b) dimension of the
decision-space; or (c) number of uncertain constraints, is small, then computing
the optimal 2-piecewise constant second stage policy can be done efficiently. One
application where the dimension of the uncertainty is large, but can be approximated by
a low-dimensional set, is weather uncertainty in air traffic flow management (see
\cite{BertsimasCaramanisTAC}).

\subsubsection{Network Design}
In Atamturk and Zhang \cite{AtamturkZhang2004}, the authors consider
two-stage robust network flow and design, where the demand vector is
uncertain. This work deals with computing the optimal second stage
adaptability, and characterizing the first-stage feasible set of
decisions. While this set is convex, solving the separation problem,
and hence optimizing over it, can be NP-hard, even for the two-stage
network flow problem.

Given a directed graph $G=(V,E)$, and a demand vector $\mb{d} \in
\R^{V}$, where the edges are partitioned into first-stage and
second-stage decisions, $E = E_1 \cup E_2$, we want to obtain an
expression for the feasible first-stage decisions. We define some
notation first. Given a set of nodes, $S \subseteq V$, let
$\delta^+(S), \delta^-(S)$, denote the set of arcs into and out of
the set $S$, respectively. Then, denote the set of flows on the
graph satisfying the demand by
$$
{\cal P}_{\mb{d}} \bydef \{ \mb{x} \in \R_+^{E} \,:\,
\mb{x}(\delta^+(i)) - \mb{x}(\delta^-(i)) \geq d_i, \,\, \forall i
\in V\}.
$$
If the demand vector $\mb{d}$ is only known to lie in a given
compact set ${\cal U} \subseteq \R^V$, then the set of flows
satisfying every possible demand vector is given by the intersection
${\cal P} = \bigcap_{\mb{d} \in {\cal U}} {\cal P}_{\mb{d}}$. If the
edge set $E$ is partitioned $E=E_1 \cup E_2$ into first and
second-stage flow variables, then the set of first-stage-feasible
vectors is:
$$
{\cal P}(E_1) \bydef \bigcap_{\mb{d} \in {\cal U}}
\mbox{Proj}_{E_1}{\cal P}_{\mb{d}},
$$
where $\mbox{Proj}_{E_1}{\cal P}_{\mb{d}} \bydef \{\mb{x}_{E_1}
\,:\, (\mb{x}_{E_1},\mb{x}_{E_2}) \in {\cal P}_{\mb{d}}\}$. Then we
have:
\begin{theorem}[\cite{AtamturkZhang2004}] A vector $\mb{x}_{E_1}$ is
an element of ${\cal P}(E_1)$ iff
$\mb{x}_{E_1}(\delta^+(S))-\mb{x}_{E_1}(\delta^-(S)) \geq \zeta_S$,
for all subsets $S \subseteq V$ such that $\delta^+(S) \subseteq
E_1$, where we have defined $\zeta_S \bydef
\max\{\mb{d}(S)\,:\,\mb{d} \in {\cal U}\}$.
\end{theorem}
The authors then show that for both the budget-restricted
uncertainty model, ${\cal U} = \{\mb{d}\,:\, \sum_{i \in V} \pi_i
d_i \leq \pi_0,\,\, \bar{\mb{d}} - \mb{h} \leq \mb{d} \leq
\bar{\mb{d}} + \mb{h}\}$, and the cardinality-restricted uncertainty
model, ${\cal U} = \{\mb{d} \,:\, \sum_{i \in V} \lceil |d_i -
\bar{d}_i| \setminus h_i \rceil \leq \Gamma, \,\, \bar{\mb{d}} -
\mb{h} \leq \mb{d} \leq \bar{\mb{d}} + \mb{h}\}$, the separation
problem for the set ${\cal P}(E_1)$ is NP-hard:
\begin{theorem}[\cite{AtamturkZhang2004}] For both classes of
uncertainty sets given above, the separation problem for ${\cal
P}(E_1)$ is NP-hard for bipartite $G(V,B)$.
\end{theorem}
These results extend also to the framework of two-stage network
design problems, where the capacities of the edges are also part of
the optimization. If the second stage network topology is totally
ordered, or an arborescence, then the separation problem becomes
tractable.

\section{Applications of Robust Optimization}\label{appsec}

In this section, we examine several applications
approached by Robust Optimization techniques.

\subsection{Portfolio optimization}\label{rpossec}

One of the central problems in finance is how to allocate monetary
resources across risky assets. This problem has received
considerable attention from the Robust Optimization community and a
wide array of models for robustness have been explored in the
literature.


\subsubsection{Uncertainty models for return mean and covariance}\label{marksssec}

The classical work of Markowitz (\cite{Markowitz52, Markowitz59})
served as the genesis for modern portfolio theory. The canonical
problem is to allocate wealth across $n$ risky assets with mean
returns $\v{\mu}\in\field{R}{n}$ and return covariance matrix
$\v{\Sigma}\in\field{S}{n}_{++}$ over a weight vector
$\v{w}\in\field{R}{n}$. Two versions of the problem arise; first,
the \emph{minimum variance problem}, i.e.,
\begin{eqnarray}\label{mark1}
\min\left\{\v{w}^{\top}\v{\Sigma}\v{w} \ : \
\v{\mu}^{\top}\v{w}\geq r, \ \v{w}\in\mathcal{W}\right\},
\end{eqnarray}
or, alternatively, the \emph{maximum return problem}, i.e.,
\begin{eqnarray}\label{mark2}
\max\left\{\v{\mu}^{\top}\v{w} \ : \
\v{w}^{\top}\v{\Sigma}\v{w}\leq \sigma^2, \v{w}\in\mathcal{W}\right\}.
\end{eqnarray}
Here, $r$ and $\sigma$ are investor-specified constants, and
$\mathcal{W}$ represents the set of acceptable weight vectors
($\mathcal{W}$ typically contains the normalization constraint
$\v{e}^{\top}\v{w}=1$ and often has ``no short-sales'' constraints,
i.e., $w_i\geq 0, \ i=1,\ldots,n$, among others).

While this framework proposed by Markowitz revolutionized the
financial world, particularly for the resulting insights in trading
off \emph{risk} (variance) and \emph{return}, a fundamental drawback
from the practitioner's perspective is that $\v{\mu}$ and
$\v{\Sigma}$ are rarely known with complete precision. In turn,
optimization algorithms tend to exacerbate this problem by finding
solutions that are ``extreme'' allocations and, in turn, very
sensitive to small perturbations in the parameter estimates.

Robust models for the mean and covariance information are a natural
way to alleviate this difficulty, and they have been explored by
numerous researchers. Lobo and Boyd \cite{LoboBoyd00} propose box,
ellipsoidal, and other uncertainty sets for $\v{\mu}$ and
$\v{\Sigma}$. For example, the box uncertainty sets have the form
\begin{eqnarray*}
\mathcal{M} &=& \left\{\v{\mu}\in\field{R}{n} \ \big| \
\underline{\mu}_i\leq\mu\leq\overline{\mu}_i, \ i=1,\ldots,n\right\}
\\
\mathcal{S} &=& \left\{\v{\Sigma}\in\field{S}{n}_+ \ | \
\underline{\Sigma}_{ij}\leq\Sigma_{ij}\leq\overline{\Sigma}_{ij}, \
i=1,\ldots,n, \ j=1,\ldots,n\right\}.
\end{eqnarray*}
In turn, with these uncertainty structures, they provide a
polynomial-time cutting plane algorithm for solving robust variants
of Problems (\ref{mark1}) and (\ref{mark2}), e.g., the \emph{robust
minimum variance problem}
\begin{eqnarray}\label{rmvp}
\min\left\{\sup_{\v{\Sigma}\in\mathcal{S}}\v{w}^{\top}\v{\Sigma}\v{w} \ : \ \inf_{\v{\mu}\in\mathcal{M}}\v{\mu}^{\top}\v{w}
\geq r, \ \v{w}\in\mathcal{W}\right\}.
\end{eqnarray}

Costa and Paiva \cite{CostaPaiva02} propose uncertainty structures
of the form $\mathcal{M}=\ch{\v{\mu}_1,\ldots,\v{\mu}_k}$,
$\mathcal{S}=\ch{\v{\Sigma}_1,\ldots,\v{\Sigma}_k},$
and formulate robust counterparts of (\ref{mark1}) and (\ref{mark2})
as optimization problems over linear matrix inequalities.

T\"{u}t\"{u}nc\"{u} and Koenig \cite{TutuncuKoening04} focus on the
case of box uncertainty sets for $\v{\mu}$ and $\v{\Sigma}$ as well
and show that Problem (\ref{rmvp}) is equivalent to the \emph{robust
risk-adjusted return problem}
\begin{eqnarray}\label{rrarp}
\max\left\{\inf\limits_{\v{\mu}\in\mathcal{M}, \
\v{\Sigma}\in\mathcal{S}}
\left\{\v{\mu}^{\top}\v{w}-\lambda\v{w}^{\top}\v{\Sigma}\v{w}\right\} \ : \ \v{w}\in\mathcal{W}\right\},
\end{eqnarray}
where $\lambda\geq 0$ is an investor-specified risk factor. They are
able to show that this is a saddle-point problem, and they use an
algorithm of Halld\'{o}rsson and T\"{u}t\"{u}nc\"{u}
\cite{HalldorssonTutuncu03} to compute robust efficient frontiers
for this portfolio problem.

\subsubsection{Distributional uncertainty models}\label{distsssec}

Less has been said by the Robust Optimization community about
\emph{distributional} uncertainty for the return vector in portfolio
optimization, perhaps due to the popularity of the classical
mean-variance framework of Markowitz. Nonetheless, some work has
been done in this regard. Some interesting research on that front is
that of El Ghaoui et al. \cite{ElGhaouiOksOustry03}, who examine the
problem of worst-case \emph{value-at-risk} (VaR) over portfolios
with risky returns belonging to a restricted class of probability
distributions. The $\epsilon$-VaR for a portfolio $\v{w}$ with risky
returns $\tilde{\v{r}}$ obeying a distribution $\mathbb{P}$ is the
optimal value of the problem
\begin{eqnarray}\label{vardef}
\min\left\{\gamma \ : \ \prob{\gamma\leq
-\tilde{\v{r}}^{\top}\v{w}}\leq\epsilon\right\}.
\end{eqnarray}
In turn, the authors in \cite{ElGhaouiOksOustry03} approach the
worst-case VaR problem, i.e.,
\begin{eqnarray}\label{worstcasevar}
\min\left\{V_{\mathcal{P}}(\v{w}) \ : \ \v{w}\in\mathcal{W}\right\},
\end{eqnarray}
where
\begin{eqnarray}
V_{\mathcal{P}}(\v{w}) := \begin{Bmatrix} \text{minimize} & \gamma \\
\text{subject to} & \sup\limits_{\mathbb{P}\in\mathcal{P}}
\prob{\gamma\leq -\tilde{\v{r}}^{\top}\v{w}}\leq\epsilon
\end{Bmatrix}.
\end{eqnarray}
In particular, the authors first focus on the distributional family
$\mathcal{P}$ with fixed mean $\v{\mu}$ and covariance
$\v{\Sigma}\succ\v{0}$. From a tight Chebyshev bound due to
Bertsimas and Popescu \cite{BertsimasPopescu04}, it was known that
(\ref{worstcasevar}) is equivalent to the SOCP
\begin{eqnarray*}
\min\left\{\gamma \ : \
\kappa(\epsilon)\|\v{\Sigma}^{1/2}\v{w}\|_2-\v{\mu}^{\top}\v{w} \leq
\gamma\right\},
\end{eqnarray*}
where $\kappa(\epsilon)=\sqrt{(1-\epsilon)/\epsilon}$; in
\cite{ElGhaouiOksOustry03}, however, the authors also show
equivalence of (\ref{worstcasevar}) to an SDP, and this allows them
to extend to the case of uncertainty in the moment information.
Specifically, when the supremum in (\ref{worstcasevar}) is taken
over all distributions with mean and covariance known only to belong
within $\mathcal{U}$, i.e., $(\v{\mu},\v{\Sigma})\in\mathcal{U}$,
\cite{ElGhaouiOksOustry03} shows the following:
\begin{enumerate}
\item When
$\mathcal{U}=\ch{(\v{\mu}_1,\v{\Sigma}_1),\ldots,(\v{\mu}_l,\v{\Sigma}_l)}$,
then (\ref{worstcasevar}) is SOCP-representable.

\item When $\mathcal{U}$ is a set of component-wise box constraints
on $\v{\mu}$ and $\v{\Sigma}$, then (\ref{worstcasevar}) is
SDP-representable.
\end{enumerate}

One interesting extension in \cite{ElGhaouiOksOustry03} is
restricting the distributional family to be sufficiently ``close''
to some reference probability distribution $\mathbb{P}_0$. In
particular, the authors show that the inclusion of an entropy
constraint
\begin{eqnarray*}
\int\log\frac{d\mathbb{P}}{d\mathbb{P}_0}d\mathbb{P} &\leq& d
\end{eqnarray*}
in (\ref{worstcasevar}) still leads to an SOCP-representable
problem, with $\kappa(\epsilon)$ modified to a new value
$\kappa(\epsilon,d)$ (for the details, see
\cite{ElGhaouiOksOustry03}). Thus, imposing this smoothness
condition on the distributional family only requires modification of
the risk factor.

Pinar and T\"{u}t\"{u}nc\"{u} \cite{PinarTutuncu05} study a
distribution-free model for near-arbitrage opportunities, which they
term \emph{robust profit opportunities}. The idea is as follows: a
portfolio $\v{w}$ on risky assets with (known) mean $\v{\mu}$ and
covariance $\v{\Sigma}$ is an arbitrage opportunity if (1)
$\v{\mu}^{\top}\v{w} \geq 0$, (2) $\v{w}^{\top}\v{\Sigma}\v{w} = 0$,
and (3) $\v{e}^{\top}\v{w} < 0$. The first condition implies an
expected positive return, the second implies a guaranteed return
(zero variance), and the final condition states that the portfolio
can be formed with a negative initial investment (loan).

In an efficient market, pure arbitrage opportunities cannot exist;
instead, the authors seek \emph{robust profit opportunities at level
$\theta$}, i.e., portfolios $\v{w}$ such that
\begin{eqnarray}\label{rprofit}
\v{\mu}^{\top}\v{w}-\theta\sqrt{\v{w}^{\top}\v{\Sigma}\v{w}} \geq
0, \ &\text{and}& \
\v{e}^{\top}\v{x} < 0.
\end{eqnarray}
The rationale for the system (\ref{rprofit}) is similar to the reasoning from Ben-Tal and Nemirovski \cite{BenTalNemirovski00} discussed earlier on approximations to chance constraints. Namely, under some assumptions on the distribution (boundedness and independence across the assets), portfolios that satisfy (\ref{rprofit}) have a positive return with probability at least $1-e^{-\theta^2/2}$. The authors in
\cite{PinarTutuncu05} then attempt to solve the
\emph{maximum-$\theta$ robust profit opportunity problem}:
\begin{eqnarray}\label{maxrprofit}
\sup\limits_{\theta,\v{w}}\left\{\theta \ : \
\v{\mu}^{\top}\v{w}-\theta\sqrt{\v{w}^{\top}\v{\Sigma}\v{w}} \geq 0, \ \v{e}^{\top}\v{w} < 0\right\}.
\end{eqnarray}
They then show that (\ref{maxrprofit}) is equivalent to a convex quadratic optimization problem and, under mild assumptions, has a closed-form solution.

%

Along this vein, Popescu \cite{Popescu07} has considered the problem of maximizing expected utility in a distributional-robust way when only the mean and covariance of the distribution are known. Specifically, \cite{Popescu07} shows that the problem
\begin{eqnarray}\label{popeq}
\min\limits_{\v{R}\sim (\v{\mu},\v{\Sigma})} \ev{\v{R}}{u(\v{x}'\v{R})},
\end{eqnarray}
where $u$ is any utility function and $\v{\mu}$ and $\v{\Sigma}$ denote the mean and covariance, respectively, of the random return $\v{R}$, reduces to a three-point problem. \cite{Popescu07} then shows how to optimize over this robust objective (\ref{popeq}) using quadratic programming.

\subsubsection{Robust factor models}\label{robfactorsssec}

A common practice in modeling market return dynamics is to use a
so-called \emph{factor model} of the form
\begin{eqnarray}\label{factormodel}
\tilde{\v{r}} &=& \v{\mu} + \v{V}^{\top}\v{f} + \v{\epsilon},
\end{eqnarray}
where $\tilde{\v{r}}\in\field{R}{n}$ is the vector of uncertain
returns, $\v{\mu}\in\field{R}{n}$ is an expected return vector,
$\v{f}\in\field{R}{m}$ is a vector of \emph{factor returns} driving
the model (these are typically major stock indices or other
fundamental economic indicators), $\v{V}\in\field{R}{m\times n}$ is
the \emph{factor loading matrix}, and $\v{\epsilon}\in\field{R}{n}$
is an uncertain vector of residual returns.

Robust versions of (\ref{factormodel}) have been considered by a few
authors. Goldfarb and Iyengar
\cite{GoldfarbIyengar03} consider a model with $\v{f}\in\mathcal{N}(\v{0},\v{F})$ and $\v{\epsilon}\in\mathcal{N}(\v{0},\v{D})$, then explicitly account for covariance uncertainty as:
\begin{itemize}
\item[$\bullet$] $\v{D}\in\mathcal{S}_d = \left\{\v{D} \ | \
\v{D}=\text{diag}(\v{d}), \
d_i\in\left[\underline{d}_i,\overline{d}_i\right]\right\}$

\item[$\bullet$] $\v{V}\in\mathcal{S}_v = \left\{\v{V}_0+\v{W} \ | \
\|\v{W}_i\|_g\leq\rho_i, \ i=1,\ldots,m\right\}$

\item[$\bullet$] $\v{\mu}\in\mathcal{S}_m =
\left\{\v{\mu}_0+\v{\varepsilon} \ | \ |\varepsilon|_i\leq\gamma_i,
\ i=1,\ldots,n\right\}$,
\end{itemize}
where $\v{W}_i=\v{We}_i$ and, for $\v{G}\succ\v{0}$,
$\|\v{w}\|_g=\sqrt{\v{w}^{\top}\v{Gw}}$. The authors then consider
various robust problems using this model, including robust versions
of the Markowitz problems (\ref{mark1}) and (\ref{mark2}), robust
Sharpe ratio problems, and robust value-at-risk problems, and show
that all of these problems with the uncertainty model above may be
formulated as SOCPs. The authors also show how to compute the
uncertainty parameters $\v{G}$, $\rho_i$, $\gamma_i$,
$\underline{d}_i$, $\overline{d}_i$, using historical return data
and multivariate regression based on a specific confidence
level $\omega$. Additionally, they show that a particular
ellipsoidal uncertainty model for the factor covariance matrix
$\v{F}$ can be included in the robust problems and the resulting
problem may still be formulated as an SOCP.

El Ghaoui et al. \cite{ElGhaouiOksOustry03} also consider the
problem of robust factor models. Here, the authors show how to
compute upper bounds on the robust worst-case VaR problem via SDP
for joint uncertainty models in $(\v{\mu},\v{V})$ (ellipsoidal and
matrix norm-bounded uncertainty models are considered).

\subsubsection{Multi-period robust models}

The robust portfolio models discussed heretofore have been for
single-stage problems, i.e., the investor chooses a \emph{single}
portfolio $\v{w}\in\field{R}{n}$ and has no future decisions. Some
efforts have been made on multi-stage problems. Ben-Tal et al. \cite{BenTalMargalitNemirovski2000} formulate the following, $L$-stage portfolio problem:
\begin{eqnarray}\label{btmnemportfolio}
\text{maximize} && \sum\limits_{i=1}^{n+1}r_i^Lx_i^L \nonumber \\
\text{subject to} && x_i^l = r_i^{l-1}x_i^{l-1}-y_i^l+z_i^l, \
i=1,\ldots,n, \ l=1,\ldots,L \nonumber \\
&& x_{n+1}^l = r_{n+1}^{l-1}x_{n+1}^{l-1} +
\sum\limits_{i=1}^n(1-\mu_i^l)y_i^l -
\sum\limits_{i=1}^n(1+\nu_i^l)z_i^l, \ l=1,\ldots,L \\
&& x_i^l,y_i^l,z_i^l\geq 0, \nonumber
\end{eqnarray}
Here, $x_i^l$ is the dollar amount invested in asset $i$
at time $l$ (asset $n+1$ is cash), $r_i^{l-1}$ is the uncertain return of asset $i$ from period $l-1$ to period $l$, $y_i^l$ ($z_i^l$) is the amount of asset $i$ to sell (buy) at the
beginning of period $l$, and $\mu_i^l$ ($\nu_i^l$) are the uncertain sell (buy) transaction costs of asset $i$ at period $l$.

Of course, (\ref{btmnemportfolio}) as stated is simply a linear
programming problem and contains no reference to the uncertainty in
the returns and the transaction costs. The authors note that one can
take a multi-stage stochastic programming approach to the problem,
but that such an approach may be quite difficult computationally. With tractability in mind, the authors propose an ellipsoidal uncertainty
set model (based on the mean of a period's return minus a safety
factor $\theta_l$ times the standard deviation of that period's
return, similar to \cite{PinarTutuncu05}) for the uncertain
parameters, and show how to solve a ``rolling horizon'' version of
the problem via SOCP.

%
%

Pinar and T\"{u}t\"{u}nc\"{u} \cite{PinarTutuncu05} explore a
two-period model for their robust profit opportunity problem. In
particular, they examine the problem
\begin{eqnarray}\label{pintut2period}
\sup\limits_{\v{x}_0} && \inf\limits_{\v{r}^1\in\mathcal{U}}
\sup\limits_{\theta,\v{x}^1}
\theta \nonumber \\
\text{subject to} && \v{e}^{\top}\v{x}^1 = (\v{r}^1)^{\top}\v{x}^0
\qquad\text{(self-financing constraint)} \\
&&
(\v{\mu}^2)^{\top}\v{x}^1-\theta\sqrt{(\v{x}^1)^{\top}\v{\Sigma}_2\v{x}^1}\geq
0
\nonumber \\
&& \v{e}^{\top}\v{x}^0<0 \nonumber,
\end{eqnarray}
where $\v{x}^i$ is the portfolio from time $i$ to time $i+1$,
$\v{r}^1$ is the uncertain return vector for period 1, and
$(\v{\mu}^2,\v{\Sigma}_2)$ is the mean and covariance of the return
for period 2. The tractability of (\ref{pintut2period}) depends
critically on $\mathcal{U}$, but \cite{PinarTutuncu05} derives a
solution to the problem when $\mathcal{U}$ is ellipsoidal.

\subsubsection{Computational results for robust
portfolios}\label{compportsssec}

Most of the studies on robust portfolio optimization are
corroborated by promising computational experiments. Here we provide
a short summary, by no means exhaustive, of some of the relevant
results in this vein.

\begin{itemize}
\item[$\bullet$] Ben-Tal et al. \cite{BenTalMargalitNemirovski2000}
provide results on a simulated market model, and show that their
robust approach greatly outperforms a stochastic programming
approach based on scenarios (the robust has a much lower observed
frequency of losses, always a lower standard deviation of returns,
and, in most cases, a higher mean return). Their robust approach
also compares favorably to a ``nominal'' approach that uses
expected values of the return vectors.

\item[$\bullet$] Goldfarb and Iyengar \cite{GoldfarbIyengar03}
perform detailed experiments on both simulated and real market data
and compare their robust models to ``classical'' Markowitz
portfolios. On the real market data, the robust portfolios did not
always outperform the classical approach, but, for high values of
the confidence parameter (i.e., larger uncertainty sets), the robust
portfolios had superior performance.

\item[$\bullet$] El Ghaoui et al. \cite{ElGhaouiOksOustry03} show
that their robust portfolios significantly outperform nominal
portfolios in terms of worst-case value-at-risk; their computations
are performed on real market data.

\item[$\bullet$] T\"{u}t\"{u}nc\"{u} and Koenig
\cite{TutuncuKoening04} compute robust ``efficient frontiers'' using
real-world market data. They find that the robust portfolios offer
significant improvement in worst-case return versus nominal
portfolios at the expense of a much smaller cost in expected return.

\item[$\bullet$] Erdo\u{g}an et al.
\cite{ErdoganGoldfarbIyengar2004} consider the problems of index
tracking and active portfolio management and provide detailed
numerical experiments on both. They find that the robust models of
Goldfarb and Iyengar \cite{GoldfarbIyengar03} can (a) track an index
(SP500) with much fewer assets than classical approaches (which has
implications from a transaction costs perspective) and (b) perform
well versus a benchmark (again, SP500) for active management.

\item[$\bullet$] Delage and Ye \cite{delageye09} consider a series of portfolio optimization experiments with market returns over a six-year horizon. They apply their method, which solves a distribution-robust problem with mean and covariance information based on samples (which they show can be formulated as an SDP) and show that this approach greatly outperforms an approach based on stochastic programming.

\item[$\bullet$] Ben-Tal et al. \cite{BentalBertsimasBrown08} apply
a robust model based on the theory of convex risk measures to a
real-world portfolio problem, and show that their approach can yield
significant improvements in downside risk protection at little
expense in total performance compared to classical methods.
\end{itemize}

As the above list is by no means exhaustive, we refer the reader to
the references therein for more work illustrating the computational
efficacy of robust portfolio models.

\subsection{Statistics, learning, and estimation}\label{statssec}

The process of using data to analyze or describe the parameters and
behavior of a system is inherently uncertain, so it is no surprise
that such problems have been approached from a Robust Optimization
perspective. Here we describe some of the prominent, related work.

\subsubsection{Robust Optimization and Regularization}\label{leastsqsssec}
Regularization has played an important role in many fields, including functional analysis, numerical computation, linear algebra, statistics, differential equations, to name but a few. Of interest are the properties of solutions to regularized problems. There have been a number of fundamental connections between regularization, and Robust Optimization.

El Ghaoui and Lebret consider the problem of robust least-squares solutions to systems of over-determined linear equations \cite{ElGhaouiLebret97}. Given an
over-determined system $\v{Ax}=\v{b}$, where
$\v{A}\in\field{R}{m\times n}$ and $\v{b}\in\field{R}{m}$, an
ordinary least-squares problem is $\min\limits_{\v{x}}\|\v{Ax}-\v{b}\|$.
In \cite{ElGhaouiLebret97}, the authors build explicit models to
account for uncertainty for the data $[\v{A} \ \v{b}]$. The
authors show that the solution to the $\ell^2$-regularized regression
problem, is in fact the solution to a robust optimization problem. In particular, the solution to
$$
\text{minimize} \| \v{A}\v{x} - \v{b}\| + \rho \sqrt{\|\v{x}\|_2^2 + 1},
$$
is also the solution to the robust problem
\begin{eqnarray*}\label{rls}
\min\limits_{\v{x}} && \max\limits_{\|\Delta\v{A} \
\Delta\v{b}\|_F\leq\rho}
\|(\v{A}+\Delta\v{A})\v{x}-(\v{b}+\Delta\v{b})\|,
\end{eqnarray*}
where $\|\cdot\|_F$ is the Frobenius norm of a matrix, i.e.,
$\|\v{A}\|_F=\sqrt{\trace{\v{A}^{\top}\v{A}}}$.

This result demonstrates that ``robustifying'' a solution gives us regularity properties. This has appeared in other contexts as well, for example see \cite{Lewis2002}. Drawing motivation from the robust control literature, the authors then consider extensions to structured matrix uncertainty sets, looking at the structured robust least-squares (SRLS) problem under linear, and fractional linear uncertainty structure.
%
%
%

In related work, Xu, Caramanis and Mannor \cite{XuCaramanisMannor-Lasso-TIT} consider $\ell^1$-regularized regression, commonly called Lasso, and show that this too is the solution to a robust optimization problem. Lasso has been studied extensively in statistics and signal processing (among other fields) due to its remarkable ability to recover sparsity. Recently this has attracted attention under the name of compressed sensing (see \cite{ChenDonohoSaunders,CandesTao2004}). In \cite{XuCaramanisMannor-Lasso-TIT}, the authors show that the solution to
$$
\text{minimize} \| \v{A}\v{x} - \v{b}\|_2 + \lambda \|\v{x}\|_1,
$$
is also the solution to the robust problem
\begin{eqnarray*}
\min\limits_{\v{x}} && \max\limits_{\|\Delta\v{A}\|_{\infty,2} \leq\rho}
\|(\v{A}+\Delta\v{A})\v{x}-\v{b}\|,
\end{eqnarray*}
where $\|\cdot\|_{\infty,2}$ is $\infty$-norm of the $2$-norm of the columns. Using this equivalence, they re-prove that Lasso is sparse using a new robust optimization-based explanation of this sparsity phenomenon, thus showing that sparsity is a consequence of robustness.

In \cite{XuCaramanisMannorSVM2009}, the authors consider robust Support Vector Machines (SVM) and show that like Lasso and Tikhonov-regularized regression, norm-regularized SVMs also have a hidden robustness property: their solutions are solutions to a (non-regularized) robust optimization problem. Using this connection, they prove statistical consistency of SVMs without relying on stability or VC-dimension arguments, as past proofs had done. Thus, this equivalence provides a concrete link between good learning properties of an algorithm and its robustness, and provides a new avenue for designing learning algorithms that are consistent and generalize well. For more on this, we refer to the book chapter on Robust Optimization and Machine Learning \cite{CaramanisMannorXuBookChapter}.

\subsubsection{Binary classification via linear discriminants}\label{binsssec}

Robust versions of binary classification problems are explored in
several papers. The basic problem setup is as follows: one has a
collection of data vectors associated with two classes, $\v{x}$ and
$\v{y}$, with elements of both classes belonging to $\field{R}{n}$.
The realized data for the two classes have empirical means and
covariances $(\v{\mu}_x,\v{\Sigma}_x)$ and
$(\v{\mu}_y,\v{\Sigma}_y)$, respectively. Based on the observed
data, we wish to find a linear decision rule for deciding, with high
probability, to which class future observations belong. In other
words, we wish to find a hyperplane
$\mathcal{H}(\v{a},b)=\left\{\v{z}\in\field{R}{n} \ | \
\v{a}^{\top}\v{z}=b\right\}$, with future classifications on new
data $\v{z}$ depending on the sign of $\v{a}^{\top}\v{z}-b$ such
that the misclassification probability is as low as possible. (We direct the interested reader to Chapter 12 of Ben-Tal et al. \cite{BentalElGhaouiNemirovskiBook} for more discussion on RO in classification problems).

Lanckriet et al. \cite{Lanckriet02} approach this problem first from
the approach of distributional robustness. In particular, they
assume the means and covariances are known exactly, but nothing else
about the distribution is known. In particular, the \emph{Minimax Probability
Machine} (MPM) finds a separating hyperplane $(\v{a},b)$ to the
problem
\begin{eqnarray}\label{mpm}
\max\left\{\alpha \ : \ \inf\limits_{\v{x}\sim(\v{\mu}_x,\v{\Sigma}_x)}
\prob{\v{a}^{\top}\v{x}\geq b} \geq\alpha, \ \inf\limits_{\v{y}\sim(\v{\mu}_y,\v{\Sigma}_y)}
\prob{\v{a}^{\top}\v{y}\leq b} \geq\alpha\right\},
\end{eqnarray}
where the notation $\v{x}\sim(\v{\mu}_x,\v{\Sigma}_x)$ means the inf
is taken with respect to all distributions with mean $\v{\mu}_x$ and
covariance $\v{\Sigma}_x$. The authors then show that (\ref{mpm})
can be solved via SOCP. The authors then go on to show that in the case when the means and covariances themselves belong to an uncertainty set defined as follows
\begin{eqnarray}
\mathcal{X} &=& \left\{(\v{\mu}_x,\v{\Sigma}_x) \ | \
(\v{\mu}_x-\v{\mu}_x^0)^{\top}\v{\Sigma}_x^{-1}(\v{\mu}_x-\v{\mu}_x^0)
\leq\nu^2, \ \|\v{\Sigma}_x-\v{\Sigma}_x^0\|_F\leq\rho\right\}, \label{xdef} \\
\mathcal{Y} &=& \left\{(\v{\mu}_y,\v{\Sigma}_y) \ | \
(\v{\mu}_y-\v{\mu}_y^0)^{\top}\v{\Sigma}_y^{-1}(\v{\mu}_y-\v{\mu}_y^0)
\leq\nu^2, \ \|\v{\Sigma}_y-\v{\Sigma}_y^0\|_F\leq\rho\right\},
\label{ydef}
\end{eqnarray}
that the problem reduces to an equivalent MPM of the form of (\ref{mpm}).

Another technique for linear classification is based on so-called
\emph{Fisher discriminant analysis} (FDA) \cite{Fisher36}. For
random variables belonging to class $\v{x}$ or class $\v{y}$,
respectively, and a separating hyperplane $\v{a}$, this approach
attempts to maximize the Fisher discriminant ratio
\begin{eqnarray}\label{fdr}
f(\v{a},\v{\mu}_x,\v{\mu}_y,\v{\Sigma}_x,\v{\Sigma}_y) &:=&
\frac{\left(\v{a}^{\top}(\v{\mu}_x-\v{\mu}_y)\right)^2}{\v{a}^{\top}\left(\v{\Sigma}_x
+ \v{\Sigma}_y\right)\v{a}},
\end{eqnarray}
where the means and covariances, as before, are denoted by
$(\v{\mu}_x,\v{\Sigma}_x)$ and $(\v{\mu}_y,\v{\Sigma}_y)$. The
Fisher discriminant ratio can be thought of as a ``signal-to-noise''
ratio for the classifier, and the discriminant
\begin{eqnarray*}
\v{a}^{\text{nom}} &:=&
\left(\v{\Sigma}_x+\v{\Sigma}_y\right)^{-1}(\v{\mu}_x - \v{\mu}_y)
\end{eqnarray*}
gives the maximum value of this ratio. Kim et al.
\cite{KimMagnaniBoyd05} consider the \emph{robust Fisher linear
discriminant problem}
\begin{eqnarray}\label{robflda}
\text{maximize}_{\v{a}\neq\v{0}} &&
\min\limits_{\left(\v{\mu}_x,\v{\mu}_y,
\v{\Sigma}_x,\v{\Sigma}_y\right)\in\mathcal{U}}
f(\v{a},\v{\mu}_x,\v{\mu}_y,\v{\Sigma}_x,\v{\Sigma}_y),
\end{eqnarray}
where $\mathcal{U}$ is a convex uncertainty set for the mean and
covariance parameters. The main result is that if $\mathcal{U}$ is a convex set, then the discriminant
\begin{eqnarray*}
\v{a}^* &:=&
\left(\v{\Sigma}_x^*+\v{\Sigma}_y^*\right)^{-1}(\v{\mu}_x^* -
\v{\mu}_y^*)
\end{eqnarray*}
is optimal to the Robust Fisher linear discriminant problem
(\ref{robflda}), where
$(\v{\mu}_x^*,\v{\mu}_y^*,\v{\Sigma}_x^*,\v{\Sigma}_y^*)$ is any
optimal solution to the convex optimization problem:
\begin{eqnarray*}
\min\left\{
(\v{\mu}_x-\v{\mu}_y)^{\top}(\v{\Sigma}_x+\v{\Sigma}_y)^{-1}(\v{\mu}_x-\v{\mu}_y)
\ : \
(\v{\mu}_x,\v{\mu}_y,\v{\Sigma}_x,\v{\Sigma}_y)\in\mathcal{U}\right\}.
\end{eqnarray*}
The result is general in the sense that no structural properties, other than convexity, are imposed on the uncertainty set $\mathcal{U}$.

Other work using robust optimization for classification and
learning, includes that of Shivaswamy, Bhattacharyya and Smola
\cite{ShivaswamyBhattacharyyaSmola2006} where they consider SOCP
approaches for handling missing and uncertain data, and also
Caramanis and Mannor \cite{CaramanisMannorCOLT2008}, where robust
optimization is used to obtain a model for uncertainty in the
label of the training data.

\subsubsection{Parameter estimation}\label{maxlikesssec}

Calafiore and El Ghaoui \cite{CalafioreElGhaoui01} consider the
problem of maximum likelihood estimation for linear models when
there is uncertainty in the underlying mean and covariance
parameters. Specifically, they consider the problem of estimating
the mean $\bar{\v{x}}$ of an unknown parameter $\v{x}$ with prior
distribution $\mathcal{N}(\bar{\v{x}},\v{P}(\v{\Delta}_p))$. In
addition, we have an observations vector
$\v{y}\sim\mathcal{N}(\bar{\v{y}},\v{D}(\v{\Delta}_d))$, independent
of $\v{x}$, where the mean satisfies the linear model
\begin{eqnarray}
\bar{\v{y}} &=& \v{C}(\v{\Delta}_c)\bar{\v{x}}.
\end{eqnarray}
Given an \emph{a priori} estimate of $\v{x}$, denoted by $\v{x}_s$,
and a realized observation $\v{y}_s$, the problem at hand is to
determine an estimate for $\bar{\v{x}}$ which maximizes the \emph{a
posteriori} probability of the event $(\v{x}_s,\v{y}_s)$. When all
of the other data in the problem are known, due to the fact that
$\v{x}$ and $\v{y}$ are independent and normally distributed, the
maximum likelihood estimate is given by
\begin{eqnarray}\label{nommlestimate}
\bar{\v{x}}_{\text{ML}}(\v{\Delta}) &=&
\arg\min\limits_{\bar{\v{x}}}
\|F(\v{\Delta})\bar{\v{x}}-g(\v{\Delta})\|^2,
\end{eqnarray}
where
\begin{eqnarray*}
\v{\Delta} &=& \begin{bmatrix} \v{\Delta}_p^{\top} &
\v{\Delta}_d^{\top} &
\v{\Delta}_c^{\top}\end{bmatrix}^{\top}, \\
F(\v{\Delta}) &=& \begin{bmatrix}
\v{D}^{-1/2}(\v{\Delta}_d)\v{C}(\v{\Delta}_c) \\
\v{P}^{-1/2}(\v{\Delta}_p)\end{bmatrix}, \\
g(\v{\Delta}) &=& \begin{bmatrix} \v{D}^{-1/2}(\v{\Delta}_d)\v{y}_s
\\ \v{P}^{-1/2}(\v{\Delta}_p)\v{x}_s\end{bmatrix}.
\end{eqnarray*}

The authors in \cite{CalafioreElGhaoui01} consider the case with
uncertainty in the underlying parameters. In particularly, they
parameterize the uncertainty as a linear-fractional (LFT) model and
consider the uncertainty set
\begin{eqnarray}\label{delta1def}
\v{\Delta}_1 &=& \left\{\v{\Delta}\in\hat{\v{\Delta}} \ \Big| \
\|\v{\Delta}\|\leq 1\right\},
\end{eqnarray}
where $\hat{\v{\Delta}}$ is a linear subspace (e.g.,
$\field{R}{p\times q}$) and the norm is the spectral (maximum
singular value) norm. The robust or \emph{worst-case maximum
likelihood} (WCML) problem, then, is
\begin{eqnarray}\label{wcml}
\text{minimize} && \max\limits_{\v{\Delta}\in\v{\Delta}_1}
\|F(\v{\Delta})\v{x}-g(\v{\Delta})\|^2.
\end{eqnarray}
One of the main results in \cite{CalafioreElGhaoui01} is that the
WCML problem (\ref{wcml}) may be solved via an SDP formulation. When
$\hat{\v{\Delta}}=\field{R}{p\times q}$, (i.e., unstructured
uncertainty) this SDP is exact; if the underlying subspace has more
structure, however, the SDP finds an upper bound on the worst-case
maximum likelihood.

Eldar et al. \cite{EldarBenTalNemirovski} consider the problem of
estimating an unknown, deterministic parameter $\v{x}$ based on an
observed signal $\v{y}$. They assume the parameter and observations
are related by a linear model
\begin{eqnarray*}
\v{y} &=& \v{Hx}+\v{w},
\end{eqnarray*}
where $\v{w}$ is a zero-mean random vector with covariance
$\v{C}_w$. The \emph{minimum mean-squared error (MSE) problem} is
\begin{eqnarray}\label{minmse}
\min\limits_{\hat{\v{x}}} \ev{}{\|\v{x}-\hat{\v{x}}\|^2}.
\end{eqnarray}
Obviously, since $\v{x}$ is unknown, this problem cannot be directly
solved. Instead, the authors assume some partial knowledge of
$\v{x}$. Specifically, they assume that the parameter obeys
\begin{eqnarray}\label{xnormbound}
\|\v{x}\|_{\v{T}} &\leq& L,
\end{eqnarray}
where $\|\v{x}\|_{\v{T}}^2=\v{x}^{\top}\v{Tx}$ for some known,
positive definite matrix $\v{T}\in\field{S}{n}$, and $L\geq 0$. The
\emph{worst-case MSE problem} then is
\begin{eqnarray}\label{worstcasemse}
\min\limits_{\hat{\v{x}}=\v{Gy}}
\max\limits_{\left\{\|\v{x}\|_{\v{T}}\leq L\right\}}
\ev{}{\|\v{x}-\hat{\v{x}}\|^2} = \min\limits_{\hat{\v{x}}=\v{Gy}}
\max\limits_{\left\{\|\v{x}\|_{\v{T}}\leq L\right\}}
\left\{\v{x}^{\top}(\v{I}-\v{GH})^{\top}(\v{I}-\v{GH})\v{x} +
\trace{\v{GC}_w\v{G}^{\top}}\right\}.
\end{eqnarray}
Notice that this problem restricts to estimators which are linear in
the observations. \cite{EldarBenTalNemirovski} then shows that
(\ref{worstcasemse}) may be solved via SDP and, moreover, when
$\v{T}$ and $\v{C}_w$ have identical eigenvectors, that the problem
admits a closed-form solution. The authors also extend this
formulation to include uncertainty in the system matrix $\v{H}$. In
particular, they show that the robust worst-case MSE problem
\begin{eqnarray}\label{robwcmse}
\min\limits_{\hat{\v{x}}=\v{Gy}}
\max\limits_{\left\{\|\v{x}\|_{\v{T}}\leq L, \
\|\delta\v{H}\|\leq\rho\right\}} \ev{}{\|\v{x}-\hat{\v{x}}\|^2},
\end{eqnarray}
where the matrix $\v{H}+\delta\v{H}$ is now used in the system model
and the matrix norm used is the spectral norm, may also be solved
via SDP.

For other work on sparsity and statistics, and sparse covariance estimation, we refer the reader to recent work in \cite{BanerjeeElGhaouiAspremont2008}, \cite{ApremontBachElGhaoui}, and \cite{CandesLiMaWright2009}.

\subsection{Supply chain management}\label{scmssec}

Bertsimas and Thiele \cite{BertsimasThiele06} consider a robust
model for inventory control. 
They use a cardinality-constrained uncertainty set, as developed in Section
\ref{rlpsssec}.
One main contribution of \cite{BertsimasThiele06} is to show that
the robust problem has an optimal policy which
is of the $(s_k,S_k)$ form, i.e., order an amount $S_k-x_k$ if
$x_k<s_k$ and order nothing otherwise, and the authors explicitly
compute $(s_k,S_k)$. Note that this implies that the robust approach
to single-station inventory control has policies which are
structurally identical to the stochastic case, with the added
advantage that probability distributions need not be assumed in the
robust case. A further benefit shown by the
authors is that tractability of the problem readily extends to
problems with capacities and over networks, and the authors in
\cite{BertsimasThiele06} characterize the optimal policies in these
cases as well.

Ben-Tal et al. \cite{BenTalGolanyNemirovskiVial2003} propose an
adaptable robust model, in particular an AARC for an inventory
control problem in which the retailer has flexible commitments with
the supplier.
This model has adaptability explicitly integrated into
it, but computed as an \emph{affine} function of the realized
demands. Thus, they use the affine adaptable framework of Section \ref{ssec:affine} This structure allows the authors in
\cite{BenTalGolanyNemirovskiVial2003} to obtain an approach which is
not only robust and adaptable, but also computationally tractable.
The model is more general than the above discussion in that it
allows the retailer to pre-specify order levels to the supplier
(commitments), but then pays a piecewise linear penalty for the
deviation of the actual orders from this initial specification. For
the sake of brevity, we refer the reader to the paper for details.

Bienstock and \"{O}zbay \cite{BienstockOzbay06} propose a robust
model for computing basestock levels in inventory control. One of
their uncertainty models, inspired by adversarial queueing theory,
is a non-convex model with ``peaks'' in demand, and they provide a
finite algorithm based on Bender's decomposition and show promising
computational results.

\subsection{Engineering}\label{engdesssec}

Robust Optimization techniques have been applied to a wide variety
of engineering problems. Many of the relevant references have already been provided in the individual sections above, in particular in Section \ref{sec:structuretractability} and subsections therein. In this section, we briefly mention some additional work in this area. For the sake of brevity, we omit most
technical details and refer the reader to the relevant papers for
more.

Some of the many papers on robust engineering design problems are
the following.

\begin{enumerate}
\item \emph{Structural design}. Ben-Tal and Nemirovski
\cite{BenTalNemTruss97} propose a robust version of a truss topology design problem in which the resulting truss structures have stable
performance across a family of loading scenarios. They derive an SDP approach to solving this robust design problem.

\item \emph{Circuit design}. Boyd et al. \cite{BoydKimPatilHorowitz05}
and Patil et al. \cite{PatilYunKimCheungHorowitzBoyd05} consider the
problem of minimizing delay in digital circuits when the underlying
gate delays are not known exactly. They show how to approach such
problems using geometric programming. See also
\cite{ManiSinghOrshansky2006,ManiCaramanisOrshansky2007,SinghHeCaramanisOrshansky2009},
already discussed above.

\item \emph{Power control in wireless channels}. Hsiung et al.
\cite{HsiungKimBoydWireless05} utilize a robust geometric
programming approach to approximate the problem of minimizing the
total power consumption subject to constraints on the outage
probability between receivers and transmitters in wireless channels
with lognormal fading. For more on applications to communication,
particularly the application of geometric programming, we refer the
reader to the monograph \cite{Chiang05}, and the review articles \cite{Luo2003,LuoYu2006}.
For applications to coordination schemes and power control in wireless channels, see \cite{YunCaramanis2008}.

\item \emph{Antenna design}. Lorenz and Boyd \cite{LorenzBoyd05}
consider the problem of building an array antenna with minimum
variance when the underlying array response is not known exactly.
Using an ellipsoidal uncertainty model, they show that this problem
is equivalent to an SOCP. Mutapcic et al. \cite{MutapcicKimBoyd06}
consider a beamforming design problem in which the weights cannot be
implemented exactly, but instead are known only to lie within a box
constraint. They show that the resulting design problem has the same
structure as the underlying, nominal beamforming problem and may, in
fact, be interpreted as a regularized version of this nominal
problem.

\item \emph{Control}. Notions of robustness have been widely
popular in control theory for several decades (see, e.g., Ba\c{s}ar
and Bernhard \cite{BasarBernhard1995}, and Zhou et al.
\cite{ZhouDoyleGlover1996}). Somewhat in contrast to this
literature, Bertsimas and Brown \cite{BertsimasBrown05a} explicitly
use recent RO techniques to develop a tractable approach to
constrained linear-quadratic control problems.

\item \emph{Simulation Based Optimization in Engineering}. In stark contrast to many of the problems we have thus-far described, many engineering design problems do not have characteristics captured by an easily-evaluated and manipulated functional form. Instead, for many problems, the physical properties of a system can often only be described by numerical simulation. In \cite{BertsimasNohadaniTeo2007}, Bertsimas, Nohadani and Teo present a framework for robust optimization in exactly this setting, and describe an application of their robust optimization method for electromagnetic scattering problems.

\end{enumerate}

\section{Future directions}\label{futsec}

The goal of this paper has been to survey the known landscape of the
theory and applications of RO. Some of the unknown questions
critical to the development of this field are the following:

\begin{enumerate}
\item\emph{Tractability of adaptable RO}. While in some very special
cases, we have known, tractable approaches to multi-stage RO, these
are still quite limited, and it is fair to say that most adaptable
RO problems currently remain intractable. The most pressing research
directions in this vein, then, relate to tractability, so that a
similarly successful theory can be developed as in single-stage
static Robust Optimization.

\item\emph{Characterizing the price of robustness}. Some work (e.g.,
\cite{BertsimasSim04a, XuMannor2006}) has explored the cost, in
terms of optimality from the nominal solution, associated with
robustness. These studies, however, have been largely empirical.
Of interest are theoretical bounds to gain an understanding of
when robustness is cheap or expensive.

\item\emph{Further developing RO from a data-driven perspective}. While some
RO approaches build uncertainty sets directly from data, most of the
models in the Robust Optimization literature are not directly
connected to data. Recent work on this issue (\cite{delageye09}, \cite{BertsimasBrown05b}) have started to lay a foundation to this perspective.
Further developing a data-driven theory of RO is interesting from a
theoretical perspective, and also compelling in a practical sense,
as many real-world applications are data-rich.
\end{enumerate}

\renewcommand{\baselinestretch}{1.00}
\small
\bibliographystyle{plain}
\bibliography{rosurvey}

\end{document}